\documentclass{amsart}

\usepackage{amssymb}
\usepackage{graphicx}

\newtheorem{theorem}{Theorem}[section]
\newtheorem{lemma}[theorem]{Lemma}
\newtheorem{proposition}[theorem]{Proposition}

\newtheorem{_definition}[theorem]{Definition}
\newenvironment{definition}{\begin{_definition}\rm}{\end{_definition}}

\newtheorem{_remark}[theorem]{\it Remark}
\newenvironment{remark}{\begin{_remark}\rm}{\end{_remark}}

\newtheorem{_claim}[theorem]{Claim}
\newenvironment{claim}{\begin{_claim}\rm}{\end{_claim}}

\numberwithin{equation}{section}
\numberwithin{table}{section}
\numberwithin{figure}{section}
\renewcommand{\qed}{\hfill {$\Box$}}

\newcommand{\A}{\mathord{\mathbb A}}
\newcommand{\C}{\mathord{\mathbb C}}

\renewcommand{\P}{\mathord{\mathbb  P}}
\newcommand{\Q}{\mathord{\mathbb  Q}}
\newcommand{\R}{\mathord{\mathbb R}}
\newcommand{\sphere}{\mathord{\mathbb S}}
\newcommand{\Z}{\mathord{\mathbb Z}}

\newcommand{\GGG}{\mathord{\mathcal G}}
\newcommand{\LLL}{\mathord{\mathcal L}}
\newcommand{\MMM}{\mathord{\mathcal M}}
\newcommand{\QQQ}{\mathord{\mathcal Q}}

\newcommand{\maprightsp}[1]{\; \smash{\mathop{\; \longrightarrow \; }\limits\sp{#1}}\; }
\newcommand{\maprightsb}[1]{\; \smash{\mathop{\; \longrightarrow \; }\limits\sb{#1}}\; }

\newcommand{\maprightsurjsp}[1]{\; \smash{\mathop{\; \surj \; }\limits\sp{#1}}\; }
\newcommand{\mapdown}{\phantom{\Big\downarrow}\hskip -8pt \downarrow}

\newcommand{\inj}{\hookrightarrow}
\newcommand{\surj}{\mathbin{\to \hskip -7pt \to}}
\newcommand{\isom}{\mathbin{\,\raise -.6pt\rlap{$\to$}\raise 3.5pt \hbox{\hskip .3pt$\mathord{\sim}$}\,}}

\newcommand{\set}[2]{\{\; {#1} \; \mid \; {#2} \;  \}}
\newcommand{\shortset}[2]{\{ {#1} \,|\, {#2}   \}}

\newcommand{\wt}{\widetilde}

\newcommand{\sprime}{\sp\prime}
\newcommand{\sprimeinv}{\sp{\prime-1}}

\newcommand{\dual}{\sp{\vee}}

\newcommand{\inv}{\sp{-1}}

\newcommand{\NS}{\mathord{\rm NS}}

\newcommand{\Ker}{\operatorname{\rm Ker}\nolimits}

\newcommand{\pr}{\mathord{\rm pr}}

\newcommand{\Pt}{\P^2}

\newcommand{\rmand}{\textrm{and}}

\newcommand{\quand}{\quad\rmand\quad}

\newcommand{\Tr}{\mathord{\rm T}}
\newcommand{\ori}[1]{\wt{#1}}
\newcommand{\oriT}{\wt{\Tr}}
\newcommand{\oriLLL}{\wt{\LLL}}
\newcommand{\oriGGG}{\wt{\GGG}}
\newcommand{\oriL}{\wt{L}}
\newcommand{\Emb}{\mathord{\rm Emb}}

\newcommand{\tf}{^{\rm tf}}
\newcommand{\cohomH}{\mathord{\rm H}}
\newcommand{\homH}{\mathord{\rm H}}
\newcommand{\Gal}{\mathord{\rm Gal}}
\newcommand{\Hom}{\mathord{\rm Hom}}
\newcommand{\GL}{\mathord{GL}}
\newcommand{\SL}{\mathord{SL}}
\newcommand{\des}{\sp{\sim}}
\newcommand{\erase}[1]{}
\newcommand{\bdr}{\partial}
\newcommand{\mon}{M}
\newcommand{\pione}{\pi_1}

\newcommand{\PStextplot}[3]{\rlap{\hskip 0pt  \raise 0pt \hbox{\hskip #1pt  \raise #2pt \hbox{#3}}}}

\begin{document}

\title[Transcendental lattices]{Zariski-van Kampen method and transcendental lattices of 
certain singular $K3$ surfaces}

\author{Ken-ichiro Arima}
\address{
Department of Mathematics,
Faculty of Science,
Hokkaido University,
Sapporo 060-0810,
JAPAN
}
\email{arima@math.sci.hokudai.ac.jp
}

\author{Ichiro~Shimada}
\address{
Department of Mathematics,
Graduate School of Science,
Hiroshima University,
1-3-1~Kagamiyama, Higashi-Hiroshima, 739-8526,  JAPAN
}
\email{shimada@math.sci.hiroshima-u.ac.jp
}

\subjclass[2000]{14J28, 14H50, 14H25}

\begin{abstract}
We present a method of Zariski-van Kampen type
for  the calculation of  the transcendental lattice of 
a complex projective  surface.
As an application, we 
calculate the transcendental lattices of complex singular $K3$ surfaces
associated with an arithmetic Zariski pair of 
maximizing sextics of type $A_{10}+A_{9}$
that are defined over $\Q(\sqrt{5})$ and are conjugate to each other 
by  the action of $\Gal (\Q(\sqrt{5})/\Q)$.
\end{abstract}

\maketitle

\section{Introduction}\label{sec:Introduction}
First we prepare some terminologies about lattices.
Let $R$ be $\Z$ or $\Z_p$,
where $p$ is a prime integer or $\infty$,
$\Z_p$ is the ring of $p$-adic integers for $p<\infty$, and 
$\Z_\infty$ is the field $\R$ of real numbers.
An \emph{$R$-lattice} is a free $R$-module $L$ of finite rank
with a non-degenerate symmetric bilinear form
$$
(\phantom{a}, \phantom{a})_L \;:\; L\times L \to R.
$$
A $\Z$-lattice is simply called a \emph{lattice}.
A lattice $L$ is called \emph{even} if 
$(v,v)_L\in 2\Z$ holds for any $v\in L$.
Two lattices $L$ and $L\sprime$ are said to be \emph{in the same genus}
if the $\Z_p$-lattices $L\otimes \Z_p$ and $L\sprime\otimes \Z_p$
are isomorphic for all $p$ (including $\infty$).
Then the  set of isomorphism classes of lattices are decomposed into a disjoint union of 
\emph{genera}.
Note that, if $L$ and $L\sprime$ are in the same genus and $L$ is even,
then $L\sprime$ is also even,
because being even is a $2$-adic property.
Let $L$ be a  lattice.
Then $L$ is canonically embedded into $L\dual:=\Hom (L, \Z)$
as a submodule of finite index, and $(\phantom{a}, \phantom{a})_L$  extends to
a symmetric bilinear form
$$
(\phantom{a}, \phantom{a})_{L\dual} \;:\; L\dual\times L\dual \to \Q.
$$
Suppose that $L$ is even.
We put 
$$
D_L:=L\dual/L,
$$
and define  a quadratic form
$q_L: D_L\to \Q/2\Z$  by
$$
q_L(\bar x):=(x, x)_{L\dual}\bmod 2\Z,
\quad \textrm{where}\quad \bar x=x+L\in D_L.
$$
The pair $(D_L, q_L)$ is called the \emph{discriminant form} of $L$.
By  the following result of Nikulin (Corollary 1.9.4 in~\cite{MR525944}),
each genus of even lattices is characterized by the signature and the
discriminant form.
\begin{proposition}
Two even  lattices are in the same genus if and only if 
they have the same signature and their  discriminant forms are isomorphic.
\end{proposition}
\par
\medskip
For a $K3$ surface $X$ defined over a field $k$,
we denote by $\NS(X)$ the \emph{N\'eron-Severi lattice} of $X\otimes \bar k$;
that is, $\NS(X)$ is the lattice of numerical equivalence classes
of divisors on $X\otimes \bar k$ with the intersection paring
$\NS(X)\times \NS(X)\to \Z$.
Following the  terminology of~\cite[\S8]{MR0284440} and~\cite{MR0441982}, 
we say that a $K3$ surface $X$ defined over a field of characteristic $0$
is \emph{singular} if the rank of $\NS(X)$ attains
the possible maximum $20$.
\par
\medskip
Let $S$ be a \emph{complex} $K3$ surface.
Then the second Betti cohomology group $\cohomH^2 (S,\Z)$ is regarded as a unimodular lattice
by the cup-product,
which is even of signature $(3, 19)$.
The  N\'eron-Severi lattice $\NS(S)$ is embedded into $\cohomH^2(S, \Z)$ primitively,
because we have $\NS(S)=\cohomH^2(S, \Z)\cap \cohomH^{1,1}(S)$.
We denote by $\Tr (S)$ the orthogonal complement 
of $\NS (S)$ in $\cohomH^2(S, \Z)$,
and call $\Tr (S)$ the \emph{transcendental lattice} of $S$.
Suppose that $S$ is singular in the sense above.
Then $\Tr (S)$ is an even positive-definite lattice of rank $2$.
The Hodge decomposition
$\Tr (S)\otimes \C=\cohomH^{2,0}(S)\oplus \cohomH^{0,2}(S)$
induces a canonical orientation on $\Tr (S)$.
We denote by $\oriT(S)$ the \emph{oriented transcendental lattice} of $S$.
\par
\medskip
We denote by $[2a, b, 2c]$ the symmetric matrix
$
\left[
\begin{array}{cc}
2a & b \\
b & 2c
\end{array}
\right],
$
and put
$$
\MMM:=\set{[2a,b,2c]}{a,b, c\in \Z, a>0, c>0, 4ac-b^2>0},
$$
on which $\GL(2,\Z)$ acts by $M\mapsto g^T M g$,
where $M\in \MMM$ and $g\in\GL(2, \Z)$.
The set of isomorphism classes of even positive-definite lattices 
of rank $2$ is equal to 
$$
\LLL:=\MMM/\GL(2, \Z),
$$
while the set of isomorphism classes of even positive-definite
\emph{oriented} lattices 
of rank $2$ is equal to 
$$
\oriLLL:=\MMM/\SL(2, \Z).
$$
For a matrix $[2a,b,2c]\in \MMM$, we denote by
$\oriL[2a,b,2c]\in \oriLLL$ and $L[2a,b,2c]\in \LLL$
the isomorphism classes  represented by $[2a,b,2c]$.
\par
\medskip
In~\cite{MR0441982},
Shioda and Inose proved the following:
\begin{theorem}\label{thm:SI}
The map $S\mapsto \oriT(S)$
induces a bijection from the set of isomorphism classes of 
complex singular $K3$ surfaces $S$ to the set $\oriLLL$.
\end{theorem}
The injectivity follows from  the  Torelli theorem 
by Piatetski-Shapiro and Shafarevich~\cite{MR0284440}.
In the proof of the surjectivity,
Shioda and Inose gave an explicit construction of
the  complex singular $K3$ surface with a given oriented transcendental lattice,
and they have proved the following:
\begin{theorem}\label{cor:SI}
Every complex singular $K3$ surface is 
defined over a number field.
\end{theorem}
Let $X$ be a singular $K3$ surface defined over a number field $F$.
We denote by $\Emb (F)$ the set of embeddings of $F$ into $\C$,
and for $\sigma\in \Emb (F)$, we denote by $X^\sigma$
the \emph{complex} singular $K3$ surface $X\otimes_{F, \sigma} \C$.
We define a map
$$
\tau_X\;:\; \Emb (F) \to \oriLLL
$$
by $\tau_X(\sigma):=\oriT(X\sp\sigma)$.
Then we have the following theorem by Sch\"utt~\cite{MR2346573},
which is a generalization of a result  that had been obtained in   \cite{tssK3}.
\begin{theorem}\label{thm:SS}
Let $\GGG_X\subset\LLL$
be the genus  of all $L\in \LLL$
such that $(D_L, q_L)$ is isomorphic to $(D_{\NS(X)}, -q_{\NS(X)})$,
and let $\oriGGG_X\subset \oriLLL$
be the pull-back of $\GGG_X$ by the natural projection 
$\oriLLL\to \LLL$.
Then the image of $\tau_X$ coincides with $\oriGGG_X$.
\end{theorem}
Therefore we obtain  a surjective map
$$
\tau_X\;:\; \Emb (F) \;\surj\; \oriGGG_X.
$$
Remark that, by the  classical theory  of Gauss~\cite{MR837656}, 
we can easily calculate the oriented genus $\oriGGG_X\subset\oriLLL$ from  
the finite quadratic form $(D_{\NS(X)}, -q_{\NS(X)})$.
\par
\medskip
Let $Y$ be a geometrically reduced and irreducible projective surface defined over a number field $K$,
and let $X\to Y\otimes_K F$ be a desingularization 
of $Y\otimes_K F$ defined over a finite extension $F$ of $K$.
Suppose that $X$ is a singular $K3$ surface.
Then we can  define a map 
$$
\tau_Y\;:\; \Emb (K)  \;\surj\; \oriGGG_X
$$
by the following:
\begin{proposition}\label{prop:tauY}
The map $\tau_X: \Emb (F) \to \oriGGG_X$ factors as
$$
\Emb (F)\; \maprightsurjsp{\rho}\;  \Emb(K)\; \maprightsp{\tau_Y}\;  \oriGGG_X, 
$$
where $\rho : \Emb (F)\surj  \Emb(K)$ is the natural restriction map $\rho(\sigma):=\sigma|K$. 
\end{proposition}
The purpose of this paper is to present a method to calculate 
the map $\tau_Y$ from a defining equation of $Y$.
\par
\medskip
More generally,
we consider the following problem.
Let $S$ be a reduced irreducible complex projective surface.
For a desingularization $S\des\to S$, 
we put
$$
\cohomH^2(S\des):=\cohomH^2(S\des,\Z)/(\textrm{the torsion part}),
$$
which is regarded as a lattice by the cup-product, 
and let $\NS(S\des)\subset \cohomH^2(S\des)$
be the sublattice  of the cohomology classes of divisors on $S\des$.
We denote by 
$$
\Tr(S\des)\subset \cohomH^2(S\des)
$$
the orthogonal complement of $\NS(S\des)$ in $\cohomH^2(S\des)$.
Then we can easily see that the isomorphism class of the lattice
$\Tr(S\des)$ does not depend on the choice of the desingularization 
$S\des\to S$,
and hence we can define the \emph{transcendental lattice $\Tr (S)$ of $S$}
to be $\Tr(S\des)$. (See Lemma 3.1 of~\cite{Shioda_K3SPnew} or Proposition~\ref{prop:indep} of this paper.)
We will give  a \emph{method of  Zariski-van Kampen type}
for the calculation of  $\Tr(S)$.  
\par
\medskip
We apply our method to maximizing sextics.
Following   Persson~\cite{MR661198, MR805337}, 
we say that 
a reduced projective plane curve $C\subset \Pt$ of degree $6$ defined over
a field $k$ of characteristic $0$ is  a \emph{maximizing sextic}
if $C\otimes \bar k$ has only simple singularities and its total Milnor number 
attains the possible maximum $19$,
where $\bar k$ is the algebraic closure of $k$.
The \emph{type} of a maximizing sextic $C$ is the $ADE$-type of the singular points of $C\otimes \bar k$.
\par
\medskip
Let $C\subset \Pt$ be a maximizing sextic defined over a number field $K$.
The double covering  $Y_C\to \Pt$ branching exactly along $C$ 
is defined over $K$.
Let $X_C\to Y_C\otimes_K F$
be the minimal resolution defined over a finite extension $F$ of $K$.
Then  
$X_C$ is a singular $K3$ surface defined over $F$.
We denote by $\oriGGG_{[C]}$
the oriented genus $\oriGGG_{X_C}$.
By Proposition~\ref{prop:tauY},
we have a surjective map
$$
\tau_{[C]}:=\tau_{Y_C}\;:\; \Emb (K)\;\surj\; \oriGGG_{[C]}.
$$
\par
\medskip
As an illustration of our Zariski-van Kampen method, 
we calculate $\tau_{[C_0]}$
for  a reducible maximizing sextic $C_0$ of type $A_{10}+A_{9}$
defined over $K=\Q(\sqrt{5})$ by the homogeneous equation
\begin{equation}\label{eq:zGH}
z \cdot (\;G(x,y,z)\;+\;\alpha \cdot H(x,y,z)\;)\;=\;0,
\end{equation}
where $\alpha^2=5$ and
\erase{
\begin{eqnarray*}
G(x, y, z)&=&-9\,{x}^{4}z-14\,{x}^{3}yz+58\,{x}^{3}{z}^{2}-48\,{x}^{2}{y}^{2}z
-64\,
{x}^{2}y{z}^{2}+10\,{x}^{2}{z}^{3}+\\&&+108\,x{y}^{3}z-20\,x{y}^{2}{z}^{2}-
44\,{y}^{5}+10\,{y}^{4}z,\\
H(x, y, z)&=&5\,{x}^{4}z+10\,{x}^{3}yz-30\,{x}^{3}{z}^{2}+30\,{x}^{2}{y}^{2}z+20\,{
x}^{2}y{z}^{2}-40\,x{y}^{3}z+20\,{y}^{5}.
\end{eqnarray*}
}
\begin{eqnarray*}
G(x, y, z)&=&-9\,{x}^{4}z-14\,{x}^{3}yz+58\,{x}^{3}{z}^{2}-48\,{x}^{2}{y}^{2}z
-64\,{x}^{2}y{z}^{2}+\\&&+10\,{x}^{2}{z}^{3}+108\,x{y}^{3}z-20\,x{y}^{2}{z}^{2}-
44\,{y}^{5}+10\,{y}^{4}z,\\
H(x, y, z)&=&5\,{x}^{4}z+10\,{x}^{3}yz-30\,{x}^{3}{z}^{2}+30\,{x}^{2}{y}^{2}z+20\,{
x}^{2}y{z}^{2}-\\&&-40\,x{y}^{3}z+20\,{y}^{5}.
\end{eqnarray*}
This equation was discovered by means of   Roczen's result~\cite{MR1175728}
(see~\S\ref{sec:equation}).
\par
\medskip
We can calculate $\NS (X_{C_0})$ by  the method of Yang~\cite{MR1387816}, and 
obtain 
$$
\oriGGG_{[C_0]}\;=\;\{\;\oriL[8,3,8], \,\oriL[2,1,28]\;\}.
$$
Let $\sigma_{\pm}$ be the embeddings of $\Q(\alpha)$ into $\C$ 
given by $\sigma_\pm (\alpha)=\pm\sqrt{5}$.
We have two surjective maps 
from $\Emb (\Q(\alpha))=\{\sigma_+, \sigma_-\}$ to $\oriGGG_{[C_0]}$.
Remark that,
since the two complex maximizing sextics $C_0\sp{\sigma_+}$  and $C_0\sp{\sigma_-}$ 
cannot be distinguished by any algebraic methods,
\emph{we have to employ some transcendental method 
to determine which surjective map is the map $\tau_{[C_0]}$}.
By the method described in \S\ref{sec:ZvK} of this paper,
we obtain the following:
\begin{proposition}\label{prop:theexample}
$$
\tau_{[C_0]}(\sigma_+)=\oriL[2,1,28],
\quad 
\tau_{[C_0]}(\sigma_-)=\oriL[8,3,8].
$$
\end{proposition}
We have shown in~\cite{AZP} and~\cite{nonhomeo} that, 
for a complex maximizing sextic $C$,
the transcendental lattice $\Tr_{[C]}:=\Tr (Y_C)$
of the double covering $Y_C\to \Pt$ branching exactly along $C$  is a topological invariant
of $(\Pt, C)$.
Thus the curves $C_0^{\sigma_+}$ and $C_0^{\sigma_-}$
form an \emph{arithmetic Zariski pair}.
(See~\cite{AZP} for the definition.)
The proof of Proposition~\ref{prop:theexample}
illustrates very explicitly how the action of the Galois group of $\Q(\alpha)$ over $\Q$
affects the topology of the embedding of $C_0$ into  $\Pt$.
\par
\medskip
The first example of arithmetic Zariski pairs was discovered by 
Artal, Carmona and Cogolludo~\cite{MR2247887} in degree $12$
by means of the \emph{braid monodromy}.
It will be an interesting problem to investigate the relation between
 the  braid monodromy of a maximizing sextic $C\subset \Pt$ and our lattice invariant $\Tr_{[C]}$.
\par
\medskip
In the study of Zariski pairs of complex plane curves,
the topological fundamental groups of the complements (or its variations like  the Alexander polynomials)
have been used to distinguish the topological types.
(See, for example, \cite{MR1257321}, \cite{MR1167373} or \cite{MR1421396}
for the oldest  example of Zariski pairs of $6$-cuspidal  sextics~\cite{MR1506719, MR1507244}.)
We can calculate the fundamental groups 
$\pione (\Pt\setminus C_0^{\sigma_+})$
and 
$\pione (\Pt\setminus C_0^{\sigma_-})$
of our example
in terms of generators and relations by the classical Zariski-van Kampen theorem.
 (See, for example,~\cite{MR1341806, MR1952329}.)
It will be an interesting problem to determine whether these two groups  are isomorphic or not.
Note that,
by the theory of algebraic fundamental groups,
their profinite completions  are isomorphic. 
\par
\medskip
The plan of this paper is as follows.
In \S\ref{sec:tauY}, we prove Proposition~\ref{prop:tauY}.
In \S\ref{sec:ZvK},
we present the Zariski-van Kampen method for the calculation of the transcendental lattice
in full generality.
In \S\ref{sec:theexample}, we  apply this method to the 
complex maximizing sextics $C_0^{\sigma_\pm}$
and prove Proposition~\ref{prop:theexample}. 
In \S\ref{sec:equation}, 
we explain how we have obtained   the equation~\eqref{eq:zGH}
of $C_0$.
\par
\medskip
Thanks are due to
the referee for his/her  comments and suggestions
on the first version of  this paper.
\section{The map $\tau_Y$}\label{sec:tauY}
We recall the proof of Theorem~\ref{thm:SS},
and prove Proposition~\ref{prop:tauY}.
The main tool is Theorem~\ref{thm:SI} due to Shioda and Inose~\cite{MR0441982}.
\begin{proof}[Proof of Theorem~\ref{thm:SS}]
It is easy to see that 
the image of $\tau_X$ is contained in $\oriGGG_X$.
(See Theorem 2 in~\cite{tssK3} or Proposition 3.5 in~\cite{nonhomeo}.)
In~\cite{tssK3} and~\cite{MR2346573},
using Shioda-Inose construction, 
we constructed a singular $K3$ surface
$X^0$ defined over a number field $F^0$
such that $\NS(X)\cong \NS(X^0)$ (and hence $\oriGGG_X=\oriGGG_{X^0}$)  holds and 
that the image of $\tau_{X^0}$ coincides with 
$\oriGGG_X$.
(See also \S4 of~\cite{nonhomeo}.)
We choose an arbitrary $\sigma\in \Emb(F)$.
Then there exists $\sigma^0\in \Emb(F^0)$
such that $\tau_{X^0} (\sigma^0)=\tau_X(\sigma)$.
Since  $(X^0)^{\sigma^0}$ and $X^\sigma$ are isomorphic over $\C$ by Theorem~\ref{thm:SI}, 
there exists a number field $M\subset \C$  containing  both $\sigma^0(F^0)$ and $\sigma (F)$
such that we have an isomorphism 
$$
X^0\otimes _{F^0, \sigma^0} M\;\cong\; X\otimes _{F, \sigma} M
$$
over $M$.
Consider the commutative diagram
$$
\begin{array}{ccccccc}
&&&\Emb (F^0)& && \\
&\raise 6pt \llap{\small$\rho_{M, \sigma^0(F^0)}$}\nearrow &&& &\searrow \raise 6pt \rlap{\small$\tau_{X^0}$}& \\
\Emb (M)&&&\longrightarrow \hskip -16pt \raise 10pt \hbox{\small\makebox(0,0){$\tau_{X^0\otimes M}=\tau_{X\otimes M}$}}  
&&& \oriGGG_X=\oriGGG_{X^0}\\
&\raise -6pt \llap{\small$\rho_{M, \sigma(F)}$}\searrow &&&&\nearrow \raise -6pt \rlap{\small$\tau_{X}$}& \\
&&&\Emb (F) &&&
\end{array},
$$
where $\rho_{M, \sigma^0(F^0)}$ and $\rho_{M, \sigma(F)}$ are 
the natural surjective restriction maps.
The surjectivity of 
$\tau_{X}$ then follows from the surjectivity of $\tau_{X^0}$.
\end{proof}
\begin{proof}[Proof of Proposition~\ref{prop:tauY}]
Let $\sigma_1$ and $\sigma_2$ be elements 
of $\Emb(F)$ such that $\sigma_1|K=\sigma_2|K$.
We put
$$
\sigma_K:=\sigma_1|K=\sigma_2|K\in \Emb(K).
$$
Then the complex surfaces $X^{\sigma_1}$ and $X^{\sigma_2}$
are desingularizations of the complex surface $Y^{\sigma_K}$.
Hence Proposition~\ref{prop:tauY} follows from 
Proposition~\ref{prop:indep} below.
\end{proof}
\begin{proposition}\label{prop:indep}
Let $S_1\des$ and $S_2\des$ be two desingularizations of 
a reduced irreducible complex projective surface $S$.
Then $\Tr(S_1\des)\cong \Tr(S_2\des)$.
If $S_1\des$ and $S_2\des$ are singular $K3$ surfaces,
then $\oriT(S_1\des)\cong \oriT(S_2\des)$.
\end{proposition}
\begin{proof}
Using  a desingularization of $S_1\des\times_S S_2\des$,
we obtain a complex smooth projective surface
$\Sigma$ with birational morphisms
$\Sigma\to S_1\des$ and $\Sigma\to S_2\des$.
Since the transcendental lattice of a complex smooth projective surface
is invariant under a blowing-up,
and  any birational morphism between  smooth projective surfaces factors  into a composite of blowing-ups,
we have $\Tr(\Sigma)\cong\Tr(S_1\des)$ and $\Tr(\Sigma)\cong\Tr(S_2\des)$.
\end{proof}
\section{Zariski-van Kampen method for transcendental lattices}\label{sec:ZvK}
%
%
%
%
%
%
For a $\Z$-module $A$,
we denote by  
$$
A\tf:=A/(\textrm{the torsion part})
$$
the maximal torsion-free quotient of $A$.
If we have a bilinear form 
$A\times A\to \Z$,
then it induces a canonical bilinear form
$A\tf\times A\tf\to \Z$.
\par
\medskip
Let $S$ be a reduced irreducible complex projective surface.
Our goal is  to calculate $\Tr(S)$.
Let $\delta:S\des\to S$ be a desingularization.
We choose a reduced curve $D$ on $S$
with the following properties:
\begin{itemize}
\item[(D1)] the classes of  irreducible components of the total transform $D^\sim\subset S\des$
of $D$ span $\NS(S\des)\otimes\Q$ over $\Q$, and 
\item[(D2)] the desingularization $\delta$ induces an isomorphism $S\des\setminus D\sp\sim \cong S\setminus D$.
\end{itemize}
We put 
$$
S^0:=S\setminus D, 
$$
and 
consider the free $\Z$-module 
$$
\homH_2(S^0):=\homH_2(S^0, \Z)\tf
$$
with the intersection paring
$$
\iota\;:\; \homH_2(S^0)\times \homH_2(S^0) \to \Z.
$$
We put 
\begin{equation}\label{eq:defI}
I(S^0):=\shortset{x\in \homH_2(S^0)}{\iota(x,y)=0\;\textrm{for any}\; y\in \homH_2(S^0)},
\end{equation}
and set 
$$
V_2(S^0):= \homH_2(S^0)/I(S^0).
$$
Then $V_2(S^0)$ is a free $\Z$-module, and the intersection paring $\iota$ induces 
a non-degenerate symmetric bilinear form
$$
\bar\iota\;:\; V_2(S^0)\times V_2(S^0) \to \Z.
$$
\begin{proposition}\label{prop:VT}
The transcendental lattice $\Tr(S)=\Tr(S\des)$ is isomorphic
to the lattice $(V_2(S^0), \bar\iota)$.
\end{proposition}
\begin{proof}
By the condition (D2),
we can regard $S^0$ as a Zariski open subset of $S\des$.
Consider the homomorphism
$$
j_* \;:\; \homH_2(S^0) \to \homH_2 (S\des):=\homH_2 (S\des,\Z)\tf
$$
induced by the inclusion $j:S^0\inj S\des$.
Under the isomorphism of lattices
$$
 \homH_2 (S\des)\cong  \cohomH^2 (S\des):=\cohomH^2 (S\des,\Z)\tf
$$
induced by the Poincar\'e duality,
the image of $j_*$ is contained in $\Tr (S\des)\subset \cohomH^2 (S\des)$
by the condition (D1) on $D$.
Using the argument in the proof of Theorem 2.6 of~\cite{AZP}
or Theorem 2.1 of~\cite{nonhomeo},
we see that the homomorphism
\begin{equation*}\label{eq:surj}
j_* \;:\; \homH_2(S^0) \to \Tr (S\des)
\end{equation*}
is surjective.
Note that we have 
\begin{equation*}\label{eq:iota}
\iota(x, y)=(j_*(x), j_*(y))_T
\end{equation*}
for any $x, y\in \homH_2(S^0)$, 
where $(\phantom{a}, \phantom{a})_T$ is the cup-product on $\Tr (S\des)$.
Since $(\phantom{a}, \phantom{a})_T$ is non-degenerate,
we conclude that   $\Ker j_*=I(S^0)$.
\end{proof}
Proposition~\ref{prop:VT} shows that, in order to obtain $\Tr(S)$,
it is enough to  calculate $\homH_2(S^0)$ and $\iota$.
Enlarging $D$ if necessary,
we have a surjective morphism
$$
\phi \;:\;S^0\to U
$$
onto a Zariski open subset $U$ of an  affine line $\A^1$
such that its general fiber is a connected  Riemann surface.
By the condition (D1) on $D$,
the general fiber  of $\phi$ is non-compact.
Let $\overline{S^0}$ be a smooth irreducible projective surface containing $S^0$
as a Zariski open subset such that
$\phi$ extends to a morphism
$$
\bar\phi\;:\; \overline{S^0}\to \P^1.
$$
Let $V_1, \dots, V_M$ and $H_1, \dots, H_N$
be the irreducible components of the boundary $\overline{S^0}\setminus S^0$,
where $V_1, \dots, V_N$ are the vertical components
(that is, $\bar\phi(V_i)$ is a point),
and $H_1, \dots, H_M$ are the horizontal components
(that is, $\bar\phi(H_j)=\P^1$).
Since  the general fiber of $\phi$ is non-compact,
we have at least one horizontal component.
We put 
$$
\A^1\setminus U=\{p_1, \dots, p_m\}.
$$
Adding to $D$ some fibers of $\phi$ and making  $U$ smaller if necessary,
we can assume the following:
\begin{enumerate}
\item the surjective morphism $\phi: S^0\to U$ has only ordinary critical points,
\item\label{etale} $\bar\phi|\cup_j H_j: \cup_j H_j\to \P^1$ is \'etale over $U$, and 
\item\label{V} $V_1\cup\dots\cup V_N=\bar\phi\inv(\infty)\cup \bar\phi\inv(p_1)\cup \dots\cup \bar\phi\inv(p_m)$,
where $\{\infty\}=\P^1\setminus\A^1$.
\end{enumerate}
Note that $\bar\phi$ has no critical points on $(\cup H_j)\cap \bar\phi\inv (U)$
by the condition (\ref{etale}).
We  denote by $c_1, \dots, c_n\in U$ the critical values of $\phi$,
and put
$$
U^\sharp :=U\setminus\{c_1, \dots, c_n\}.
$$
By the assumptions,
$\phi$ is locally trivial (in the category of topological spaces and continuous maps)
over $U\sp\sharp$
with the fiber being a connected Riemann surface of genus $g$ with $r$ punctured points,
where $r>0$ is the degree of $\bar\phi|\cup_j H_j: \cup_j H_j\to \P^1$.
We then choose a base point $b\in U\sp\sharp$,
and put
$$
F_b:=\phi\inv(b).
$$
For each $p_i\in \A^1\setminus U$,
we choose a loop
$$
\lambda_i \;:\; (I, \bdr I)\to (U\sp\sharp, b)
$$
that is sufficiently smooth  and injective in the sense that 
$\lambda_i(t)=\lambda_i(t\sprime)$ holds only when $t=t\sprime$ or $\{t, t\sprime\}=\bdr I$,
and that defines the same element in $\pione(U\sp\sharp, b)$
as a simple loop (a lasso) around $p_i$ in $U\sp\sharp$.
For each critical value $c_j\in U\setminus U\sp\sharp$,
we choose a sufficiently smooth and injective  path
$$
\gamma_j \;:\; I \to U
$$
such that $\gamma_j(0)=b$, $\gamma_j(1)=c_j$ and $\gamma_j(t)\in U\sp\sharp$ for $t<1$.
We choose these loops $\lambda_i$ and paths $\gamma_j$ in such a way that
any two of them intersect only at $b$.
Then,
by a suitable self-homeomorphism of $\A^1$,
the objects $b$, $p_i$, $c_j$, $\lambda_i$ and $\gamma_j$ on $\A^1$
are mapped as in Figure~\ref{fig:A}.
%
\renewcommand{\PStextplot}[3]{\rlap{\hskip -250.000000 pt \hbox{\hskip #1pt  \raise #2pt \hbox{#3}}}}%
\begin{figure}
 \begin{center}
 \includegraphics{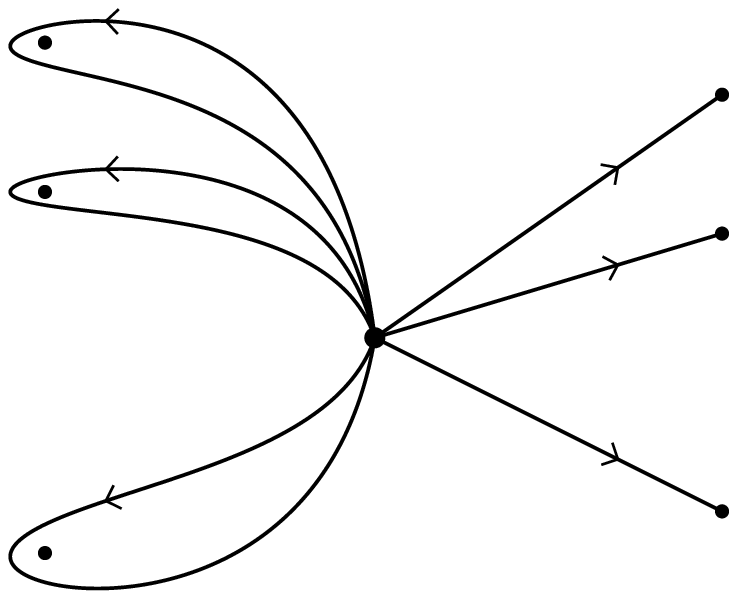}%
\PStextplot{153.000000}{85.000000}{$b$}%
\PStextplot{256.200000}{170.000000}{$c_1$}%
\PStextplot{256.200000}{130.000000}{$c_2$}%
\PStextplot{256.200000}{50.000000}{$c_n$}%
\PStextplot{61.000000}{185.000000}{$p_1$}%
\PStextplot{61.000000}{142.000000}{$p_2$}%
\PStextplot{61.000000}{38.000000}{$p_m$}%
\PStextplot{70.552000}{199.190400}{$\lambda_1$}%
\PStextplot{70.552000}{156.535200}{$\lambda_2$}%
\PStextplot{70.552000}{62.027200}{$\lambda_m$}%
\PStextplot{214.000000}{159.000000}{$\gamma_1$}%
\PStextplot{214.000000}{131.000000}{$\gamma_2$}%
\PStextplot{214.000000}{75.600000}{$\gamma_n$}%
\PStextplot{230.000000}{90.000000}{$\vdots$}%
\PStextplot{70.000000}{90.000000}{$\vdots$}%
\end{center}
\caption{}
\label{fig:A}
\end{figure}
In particular,
the union  $B$ 
of $\lambda_i$ and  $\gamma_j$
is a strong deformation retract of $U$.
Note that  $\phi$ is locally trivial over $U\setminus B$.
\par
\medskip
Let $\sphere^1$ be an oriented one-dimensional sphere.
We fix a system of oriented simple closed curves
$$
a_\nu \;:\;\sphere^1\inj  F_b \quad(\nu=1, \dots, 2g+r-1)
$$
on $F_b$ in such a way that their union
$\bigcup a_{\nu}(\sphere^1)$
is a strong deformation retract of $F_b$.
In particular, we have
$$
\homH_1 (F_b, \Z)=\bigoplus\; \Z[a_\nu],
$$
where $[a_\nu]$ is the homology class of $a_\nu$.
For each $p_i\in\A^1\setminus U$ and $a_\nu$, let
$$
\Lambda_{i, \nu}\;:\; \sphere^1\times I\inj S^0
$$
be an embedding such that the diagram
$$
\begin{array}{ccc}
\sphere^1\times I & \maprightsp{\Lambda_{i, \nu}} & S^0 \\
\llap{\scriptsize $\pr$}\mapdown & & \mapdown \rlap{\scriptsize $\phi$}\\
I &\maprightsb{\lambda_i} & U
\end{array}
$$
commutes and 
that 
$$
\Lambda_{i, \nu}| \sphere^1\times \{0\} \;:\; \sphere^1\inj F_b
$$
is equal to $a_\nu$.
We put 
$$
\mon_i(a_\nu):=\Lambda_{i, \nu}| \sphere^1\times \{1\} \;:\; \sphere^1\inj F_b,
$$
 where $\mon_i$ stands for   the \emph{monodromy} along $\lambda_i$, and 
 denote the homology class of $\mon_i(a_\nu)$ by
$$
\mon_i([a_\nu])\in \homH_1 (F_b, \Z).
$$

Let $\Theta$ be the topological space obtained from $\sphere^1\times I$
by contracting $ \sphere^1\times \{1\}$ to a point $v\in \Theta$;
that is, $\Theta$ is a cone over $\sphere^1$ with the vertex $v$.
Let $\pr: \Theta\to I$ be the natural projection.
Let  $c_j\in U\setminus U\sp\sharp$ be a critical value of $\phi$,
and let $\tilde{c}_j^1, \dots, \tilde{c}_j^m$ be the critical points of $\phi$ over $c_j$.
For each critical point $\tilde{c}_j^k\in \phi\inv(c_j)$,
we choose a \emph{thimble}
$$
\Gamma_j^k \;:\; \Theta\inj S^0
$$
along the path $\gamma_j$ corresponding to the ordinary node $\tilde{c}_j^k$ of $\phi\inv(c_j)$. 
 Namely,  the thimble $\Gamma_j^k$ is an embedding
such that 
$$
\begin{array}{ccc}
\Theta & \maprightsp{\Gamma_j^k} & S^0 \\
\llap{\scriptsize $\pr$}\mapdown & & \mapdown \rlap{\scriptsize $\phi$}\\
I &\maprightsb{\gamma_j} & U
\end{array}
$$
commutes, 
and that $\Gamma_j^k(v)=\tilde{c}_j^k$.
(See~\cite{MR592569} for thimbles and vanishing cycles.)
Then the simple closed curve
$$
\sigma_j^k:=\Gamma_j^k|\pr\inv(0)=-\bdr \Gamma_j^k\;:\; \sphere^1\inj F_b
$$
on $F_b$
represents the \emph{vanishing cycle} for the critical point $\tilde{c}_j^k$
along $\gamma_j$.
We denote its homology class by
$$
[\sigma_j^k]\in \homH_1 (F_b, \Z).
$$
We can assume that
$\Gamma_j^k$ and $\Gamma_j^{k\sprime}$ are disjoint if $k\ne k\sprime$.
\begin{remark}
There are two choices of the orientation of the thimble $\Gamma_j^k$ (and  hence
of the vanishing cycle $\sigma_j^k=-\bdr \Gamma_j^k$).
\end{remark}
\begin{figure}
\includegraphics{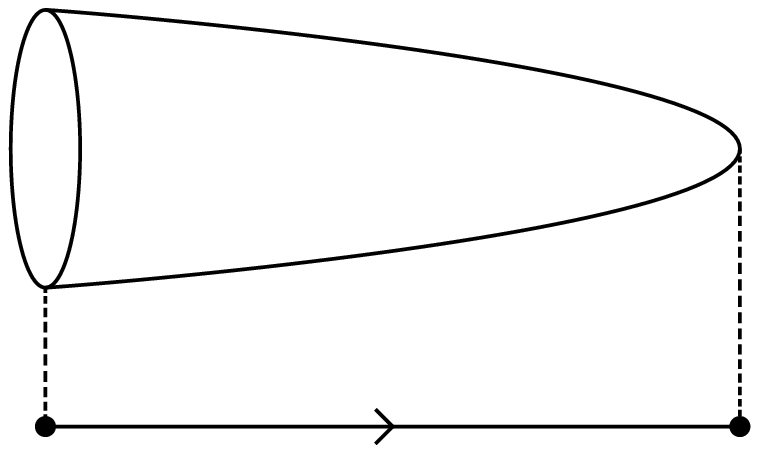}%
\rlap{\hskip -223pt \raise 3pt \hbox{$b$}}%
\rlap{\hskip -23pt \raise 3pt \hbox{$c_j$}}%
\rlap{\hskip -127pt \raise 3pt \hbox{$\gamma_j$}}%
\rlap{\hskip -14pt \raise 100pt \hbox{$\tilde{c}_j^k$}}%
\rlap{\hskip -245pt \raise 100pt \hbox{$\sigma_j^k$}}%
\rlap{\hskip -127pt \raise 100pt \hbox{$\Gamma_j^k$}}%
\caption{Thimble}
\end{figure}
Then the union
$$
F_b \;\cup\;  \bigcup \Lambda_{i, \nu} (\sphere^1\times I) \;\cup\;  \bigcup \Gamma_j^k(\Theta) 
$$
is homotopically equivalent to  $S^0$.
Since  the  $1$-dimensional CW-complex
$\bigcup a_{\nu}(\sphere^1)$
is a strong deformation retract of $F_b$, 
the homology group  $\homH_2(S^0, \Z)$ is equal to the kernel of the homomorphism
$$
\bdr\;:\; \bigoplus \Z[\Lambda_{i, \nu} ] \oplus \bigoplus \Z[\Gamma_j^k]\;\to\; \homH_1 (F_b, \Z)=\bigoplus\Z[a_\nu]
$$
given by
$$
\bdr[\Lambda_{i, \nu} ]=\mon_i([a_\nu])-[a_\nu] \quand
\bdr[\Gamma_j^k]=-[\sigma_j^k].
$$
\par
\medskip
The intersection pairing  on $\homH_2(S^0, \Z)$ is calculated by perturbing
the system  $(\lambda_i, \gamma_j)$ of loops and paths
with the base point $b$ to a system  $(\lambda\sprime_i, \gamma\sprime_j)$
with the base point $b\sprime\ne b$.
We make the perturbation in such a way that the following hold.
\begin{itemize}
\item There exists a small open disk  
$\Delta\subset U\sp\sharp$ 
containing both $b$ and $b\sprime$ such that
\begin{eqnarray*}
&&\lambda_i\inv(\Delta)=[0, s_i)\cup(1-r_i, 1], \quad
\lambda\sprimeinv_i(\Delta)=[0, s\sprime_i)\cup(1-r\sprime_i, 1], \\
&&\gamma_j\inv(\Delta)=[0, u_j), \quad\quad
\gamma\sprimeinv_j(\Delta)=[0, u\sprime_j),
\end{eqnarray*}
where $s_i, r_i, s\sprime_i, r\sprime_i, u_j, u\sprime_j$ are small positive real numbers.
\item
If $\lambda_i$ intersects $\lambda\sprime_{i\sprime}$ or $\gamma\sprime_{j\sprime}$,
then their intersection points are contained in $\Delta\setminus\{b, b\sprime\}$,
and the intersections are transverse.
\item
If $\gamma_j$ intersects $\lambda\sprime_{i\sprime}$ or $\gamma\sprime_{j\sprime}$ with $j\sprime\ne j$,
then their intersection points are contained in $\Delta\setminus\{b, b\sprime\}$,
and the intersections are transverse.
\item
Any  intersection point of $\gamma_j$ and $\gamma\sprime_j$ is either  the common end-point $c_j$, 
or a transversal intersection point contained in $\Delta\setminus\{b, b\sprime\}$.
\end{itemize}
We then perturb the topological $2$-chains
$\Lambda_{i, \nu}$ over $\lambda_{i}$  and $\Gamma_j^k$ over  $\gamma_j$
to topological $2$-chains 
$\Lambda\sprime_{i, \nu}$over $\lambda\sprime_{i}$  and $\Gamma\sp{\prime k}_j$  over $\gamma\sprime_j$,
respectively.
Let $T$ be one of $\Lambda_{i, \nu}$ or $\Gamma_j^k$,
and let $t$ be the loop or the path over which $T$ locates.
Let $T\sprime$ be one of $\Lambda\sprime_{i, \nu}$ or $\Gamma\sp{\prime k}_j$
over the loop or the path $t\sprime$.
We can make the perturbation in such a way that
$T$ and $T\sprime$ intersect transversely at each intersection points.
Suppose that
$$
t(I)\cap t\sprime(I)\cap \Delta =\{q_1, \dots, q_l\}.
$$
Then $T\cap T\sprime$ are contained in the union $\bigcup_{\mu=1}^l \phi\inv(q_\mu)$ of fibers
except for the case where $T=\Gamma_j^k$ and $T\sprime=\Gamma\sp{\prime k}_j$
for some $j$ and $k$.
If  $T=\Gamma_j^k$ and $T\sprime=\Gamma\sp{\prime k}_j$,
then $T$ and $T\sprime$ also intersect at the critical point $\tilde{c}_j^k$
transversely with the local intersection number $-1$.
(See Lemma 4.1 of~\cite{MR1433320}.)
For each $q_\mu$,
let $\theta_\mu$ and $\theta\sprime_\mu$
be the $1$-cycles on the open Riemann surface $\phi\inv(q_\mu)$ given by
\begin{eqnarray*}
\theta_\mu&:=&T\phantom{\sprime}|\sphere^1\times\{w_\mu\}\;:\; \sphere^1\to \phi\inv(q_\mu), \quad
\textrm{where $t(w_\mu)\,\,=q_\mu$},\\
\theta\sprime_\mu&:=&T\sprime|\sphere^1\times\{w\sprime_\mu\}\;:\; \sphere^1\to \phi\inv(q_\mu), \quad
\textrm{where $t\sprime(w\sprime_\mu)=q_\mu$}.
\end{eqnarray*}
We denote by $(t, t\sprime)_{\mu}$ the local intersection number of the $1$-chains $t$ and $t\sprime$
on $U$
at $q_\mu$, which is $1$ or $-1$ by the assumption on the perturbation.
We also denote by $(\theta_\mu, \theta\sprime_\mu)_{\mu}$ the intersection number
of $\theta_\mu$ and $\theta\sprime_\mu$ on the  Riemann surface $\phi\inv(q_\mu)$.
Then the intersection number $(T, T\sprime)$ of $T$ and $T\sprime$
is equal to
\begin{equation}\label{eq:TT}
(T, T\sprime)=- \sum_{\mu=1}^l  (t, t\sprime)_{\mu} (\theta_\mu, \theta\sprime_\mu)_{\mu}+\delta,
\end{equation}
where 
$$
\delta:=\begin{cases} -1 & \textrm{if $T=\Gamma_j^k$ and $T\sprime=\Gamma\sp{\prime k}_j$ for some $j$ and $k$}, \\
0 & \textrm{otherwise}.
\end{cases}
$$
The number $(\theta_\mu, \theta\sprime_\mu)_{\mu}$
is calculated as follows.
Let
$$
(\phantom{a},\phantom{a})_F\;:\; \homH_1(F_b, \Z)\times \homH_1(F_b, \Z) \to \Z
$$
be the intersection pairing  on $\homH_1(F_b, \Z)$,
which is anti-symmetric.
If $T=\Gamma_j^k$,
then the $1$-cycle $\theta_\mu$ on $\phi\inv (q_\mu)$
can be deformed to the vanishing cycle $\sigma_j^k=-\bdr \Gamma_j^k$ on $F_b$ along the  path $t|[0, w_\mu]$
in $\Delta$.
We put 
$$
[\tilde\theta_\mu]:=[\sigma_j^k]\in \homH_1(F_b, \Z).
$$
Suppose that $T=\Lambda_{i, \nu}$.
Then we have
$t\inv (\Delta)=\lambda_i\inv(\Delta)=[0, s)\cup(1-r, 1]$,
where $s$ and $r$ are small positive real numbers.
If the number $w_\mu$ such that  $t(w_\mu)=q_\mu$ is contained in $[0, s)$,
then $\theta_\mu$ can be deformed  to the $1$-cycle
$a_\nu$ on $F_b$ along the  path $t|[0, w_\mu]$
in $\Delta$.
If  $w_\mu\in (1-r, 1]$,
then $\theta_\mu$ can be deformed  to the $1$-cycle
$\mon_i(a_\nu)$ on $F_b$ along the  path $t|[w_\mu, 1]$
in $\Delta$.
We define $[\tilde\theta_\mu]\in \homH_1(F_b, \Z)$ by 
$$
[\tilde\theta_\mu]:=
\begin{cases}
[a_\nu] &\textrm{if $w_\mu\in [0, s)$}, \\
\mon_i([a_\nu]) &\textrm{if $w_\mu\in (1-r, 1]$}.
\end{cases}
$$
We define $[\tilde\theta\sprime_\mu]\in \homH_1(F_b, \Z)$
from $T\sprime$ in the same way.
Since $\phi$ is topologically  trivial over $\Delta$, 
we have 
\begin{equation}\label{eq:tt}
(\theta_\mu, \theta\sprime_\mu)_{\mu}=([\tilde\theta_\mu], [\tilde\theta\sprime_\mu])_F.
\end{equation}
The formulae~\eqref{eq:TT} and~\eqref{eq:tt}
give the intersection number $(T, T\sprime)$ of topological $2$-chains $T$ and $T\sprime$.
Even though the number $(T, T\sprime)$ depends on the choice of the perturbation,
it gives the symmetric intersection paring on $\Ker\bdr=\homH_2(S^0, \Z)$.
Thus we obtain $\homH_2(S^0)$ and $\iota$.
\section{Maximizing sextics of type $A_{10}+A_{9}$}\label{sec:theexample}
Recall from Introduction that  $\LLL$ (resp.~$\oriLLL$)
is  the set 
of isomorphism classes of even positive-definite lattices (resp.~oriented  lattices) of rank $2$.
\begin{definition}\label{def:real}
Let $\varphi: \oriLLL\to\LLL$ be the map
of forgetting orientation.
We say that $T\in \LLL$ is \emph{real} if  $\varphi\inv(T)$
consists of a single element, and
 that $\ori{T}\in \oriLLL$ is \emph{real} if 
$\varphi(\ori{T}) \in \LLL$ is real.
\end{definition}
Let $S$ be a complex singular $K3$ surface,
and let $\overline{S}$ denote $S\otimes_{\C, \bar{\phantom{z}}}\C$,
where $\bar{\phantom{z}}:\C\to\C$ is the  conjugate over $\R$.
Then $\oriT (\overline{S})$ is the reverse of $\oriT(S)$;
that is, $\varphi\inv(\varphi(\oriT(S)))=\{\oriT(S), \oriT (\overline{S})\}$.
Therefore $\oriT(S)$ is real if and only if $S$ and  $\overline{S}$ are isomorphic.
In particular, if $S$ is defined over $\R$,
then $\oriT(S)$ is real.
\begin{remark}
It is known that every element $\ori{T}$  of $\oriLLL$ is represented by a unique matrix $[2a, b, 2c]\in \MMM$
with
$$
-a<b\leq a\leq c,\;\;\textrm{with $b\geq 0$ if $a = c$}, 
$$
and  $\ori{T}$ is \emph{not} real if and only if $0<|b|<a<c$ holds.
See~\cite[Chapter 15]{MR1662447}.
\end{remark}
\par
\medskip
By the method of Yang~\cite{MR1387816} and Degtyarev~\cite{MR2357681},
we see the following facts.
(See also~\cite{AZP}.)
There are four  connected components in the moduli space of complex maximizing
sextics of type 
$A_{10}+A_{9}$.
The members of two of them are  irreducible sextics,
and their oriented transcendental lattices
are 
$$
\oriL[10, 0, 22]
\quand
\oriL[2, 0, 110]\qquad \textrm{(both are real)}.
$$
The members of  the other two are  reducible.
Each of them is a union of a line and an irreducible quintic, and 
their oriented transcendental lattices
are
$$
\oriL[8, 3, 8]
\quand
\oriL[2, 1, 28]\qquad \textrm{(both are real)}.
$$
We will consider these reducible sextics $C_0$,
whose defining equation is given by~\eqref{eq:zGH}.
\par
\begin{figure}
\hskip -160pt 
\hbox{\includegraphics[scale=.8]{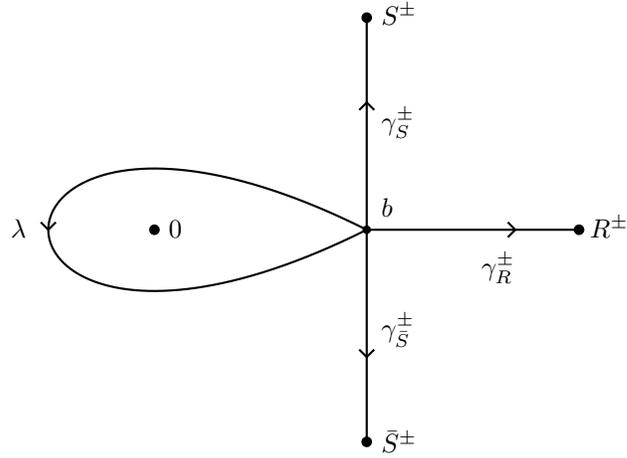}%
\rlap{\hskip -195.3pt \raise 109.5 pt \hbox{$0$}}%
\rlap{\hskip -115pt \raise 117.5 pt \hbox{$b$}}%
\rlap{\hskip -36pt \raise 109.4 pt \hbox{$R^\pm$}}%
\rlap{\hskip -115pt \raise 190 pt \hbox{$S^\pm$}}%
\rlap{\hskip -115pt \raise 28 pt \hbox{$\bar S^\pm$}}%
\rlap{\hskip -115pt \raise 150 pt \hbox{$\gamma^\pm_S$}}%
\rlap{\hskip -115pt \raise 72 pt \hbox{$\gamma^\pm_{\bar S}$}}%
\rlap{\hskip -77pt \raise 95 pt \hbox{$\gamma^\pm_R$}}%
\rlap{\hskip -255pt \raise 109.50 pt \hbox{$\lambda$}}%
}
\caption{The loop $\lambda$ and the paths $\gamma^\pm_{R}$, $\gamma^\pm_{S}$, ${\gamma}_{\bar S}^\pm$}
\label{fig1}
\end{figure}
\begin{figure}
\begin{center}
\parbox{210pt}{
\includegraphics[scale=.7]{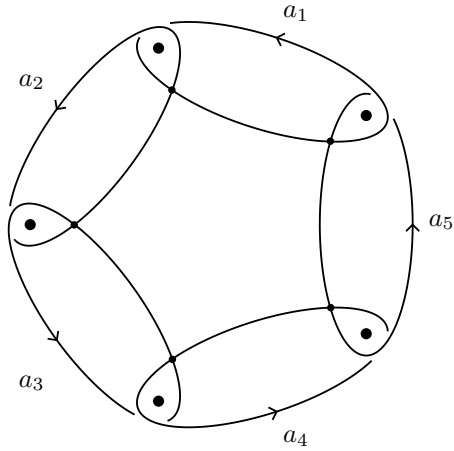}%
\rlap{\hskip -25pt \raise 105 pt \hbox{$a_5$}}%
\rlap{\hskip -80pt \raise 184 pt \hbox{$a_1$}}%
\rlap{\hskip -180pt \raise 158 pt \hbox{$a_2$}}%
\rlap{\hskip -180pt \raise 44 pt \hbox{$a_3$}}%
\rlap{\hskip -80pt \raise 23 pt \hbox{$a_4$}}%
}
\end{center}
\caption{The system of simple closed curves on $F_b$}
\label{fig2}
\end{figure}

\medskip
For simplicity, we write 
$C\sp\pm$ for $C_0^{\sigma_\pm}$,
$Y\sp\pm$ for   $Y_{C_0^{\sigma_\pm}}$ and $X^\pm$ for $X_{C_0^{\sigma_\pm}}$.
Let $D\sp\pm \subset Y\sp\pm$ be the pull-back of the union of
the lines 
$$
x=0\quand z=0
$$
on $\Pt$.
Since the singular points 
$$
[0:0:1]\;\;  \textrm{($A_{10}$)}\quand[1:0:0]\;\;  \textrm{($A_{9}$)}
$$
of $C^\pm$ are on the union of these two lines,
the curve  $D\sp\pm$ satisfies the conditions~(D1) and~(D2)
in \S\ref{sec:ZvK} for  $S=Y\sp\pm$ and $S\des=X\sp\pm$.
We denote by $W\sp\pm$ the complement of $D\sp\pm$ in $Y\sp\pm$.
(We have  denoted $W\sp\pm$  by $S^0$ in \S\ref{sec:ZvK}.)
Let $\A^2_{(y, z)}$ be the affine part of $\Pt$ given by $x\ne 0$
with the affine coordinates $(y, z)$ obtained from $[x:y:z]$ by putting $x=1$,
and let $L\subset \A^2_{(y, z)}$ be the affine line 
defined by $z=0$.
Then $W\sp\pm$ is the double cover of $\A^2_{(y, z)}\setminus L$ 
branching exactly along the union of $L$ and the smooth affine quintic curve $Q^\pm\subset \A^2_{(y, z)}$
defined by 
$$
f^{\pm}(y, z):=G(1, y, z)\pm\sqrt{5}\cdot H(1, y, z)=0.
$$
Note that $Q^\pm$ intersects $L$ only at the origin,
and   the intersection multiplicity is $5$.
Let 
$$
\pi^\pm \;:\; W^\pm \to \A_{(y, z)}^2\setminus L
$$
be the double covering.
We consider the projection
$$
p\;:\; \A^2_{(y, z)}\to \A^1_z
$$
defined by $p(y, z):=z$
onto an affine line with an affine coordinate $z$,
and the composite
$$ 
q^\pm \;:\; W^\pm \to \A_{(y, z)}^2\setminus L\to U:=\A^1_z\setminus\{0\}
$$
of $\pi^\pm$ and $p$,
which  serves as the surjective morphism $\phi$ in \S\ref{sec:ZvK}.
Calculating the discriminant  of $f^{\pm}(y, z)$ with respect to $y$, we see that
there are four  critical points of the finite covering 
$$
p|Q^\pm \;:\; Q^\pm \to \A^1_z
$$
of degree $5$.
Three of them $R^\pm, S^\pm, \overline{S}^\pm$ are simple critical values, 
where
\begin{eqnarray*}
&& R^+=0.42193...,\quad 
\;\;S^+=0.23780...+0.24431...\cdot \sqrt{-1}, \quand\\
&&R^-=0.12593...,\;\;\quad
S^-=27.542...+45.819...\cdot \sqrt{-1}.\;\;
\end{eqnarray*}
The value  $\overline{S}^\pm $  is the complex conjugate of $S^\pm$.
The critical point over $0\in \A^1_z$ is of multiplicity $5$.
The critical values of $q^\pm : W^\pm \to U$ are therefore
$R^\pm, S^\pm, \overline{S}^\pm$, and 
the fiber of $q^\pm$ over each of them has only one ordinary node.
We choose a sufficiently small positive real number $b$ as a base point on $U$,
and define  the loop $\lambda$ and the paths 
$\gamma^\pm_R$, $\gamma^\pm_S$, $\gamma^\pm _{\bar{S}}$ on $U$ as in 
Figure~\ref{fig1}.
For $z\in U$, 
we put
$$
\QQQ^\pm (z):=(p|Q^\pm)\inv (z)=p\inv (z) \cap Q^\pm, 
$$
and investigate the movement of the points $\QQQ^\pm (z)$ when $z$ moves on $U$
along the loop $\lambda$ and the paths $\gamma^\pm_R$, $\gamma^\pm_S$, $\gamma^\pm_{\bar{S}}$.
We put
$$
\A^1_y:=p\inv (b), \quad
F^\pm:=q^{\pm-1} (b)=\pi^{\pm-1}(\A^1_y)\subset W^\pm.
$$
Note that the morphism
$$
\pi^\pm|F^\pm\;:\; F^\pm\to \A^1_y
$$
 is the double covering 
branching exactly at the five points
$\QQQ^\pm (b)\subset \A^1_y$.
These branching points $\QQQ^\pm (b)$ are depicted as  big dots in Figure~\ref{fig2}.
Hence $F^\pm$ is a Riemann surface of genus $2$  minus one point.
We choose a system of oriented simple closed curves 
$$
a_\nu \;:\;\sphere^1 \inj F^\pm\qquad(\nu=1, \dots, 5)
$$
in such a way
that their images by the double covering 
$\pi^\pm|F^\pm: F^\pm\to \A^1_y$
 are  given in Figure~\ref{fig2}, 
and that  their  intersection numbers 
on $F^\pm$ are
equal to  
$$
([a_\nu],  [a_{\nu+1]})_F=-([a_{\nu+1}], [a_\nu] )_F=1
$$ 
for $\nu=1, \dots, 5$, where $a_6:=a_1$. 
(Note that $([a_\nu], [a_{\nu\sprime}])_F=0$ 
except for the case  where $|\nu-\nu\sprime|= 1$ or $\{\nu, \nu\sprime\}=\{1, 5\}$.)
Then $a_1\cup\dots\cup a_4$
is a strong deformation retract of $F^\pm$, and 
$[a_1], \dots, [a_4]$ form  a basis of
$\homH_1(F^\pm, \Z)$. 
Moreover   
 we have
$$
[a_5]=-[a_1]-[a_2]-[a_3]-[a_4].
$$
\par
\begin{figure}
\parbox{170pt}{
\includegraphics[scale=.55]{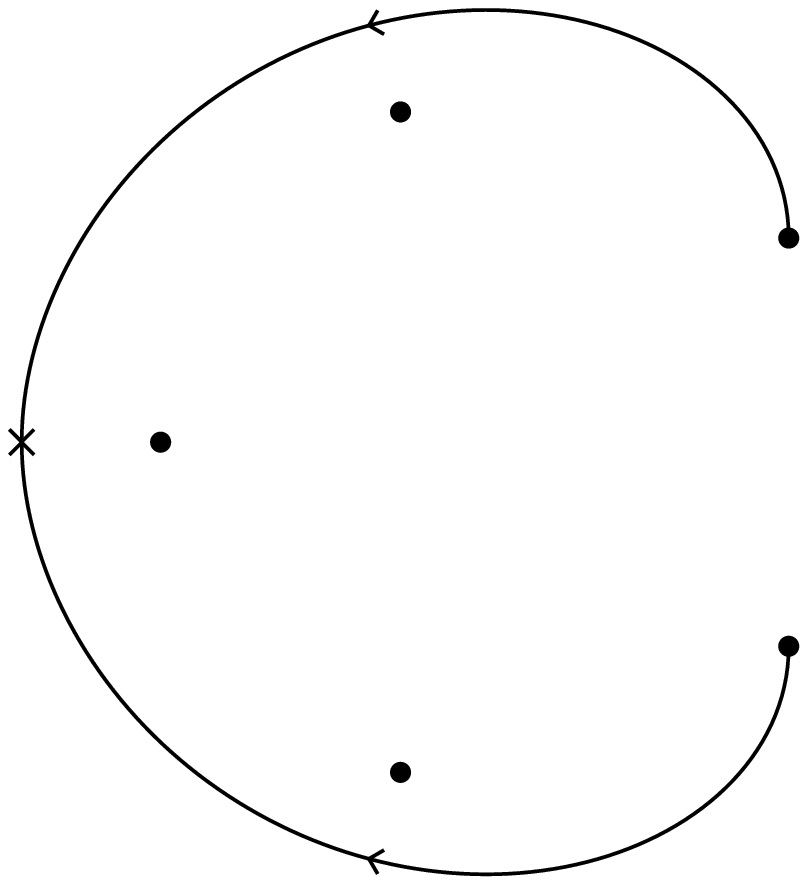}%
\\
\hfil $\QQQ^+(z)$ \hfil
}
\vrule height 100pt 
\parbox{170pt}{
\includegraphics[scale=.55]{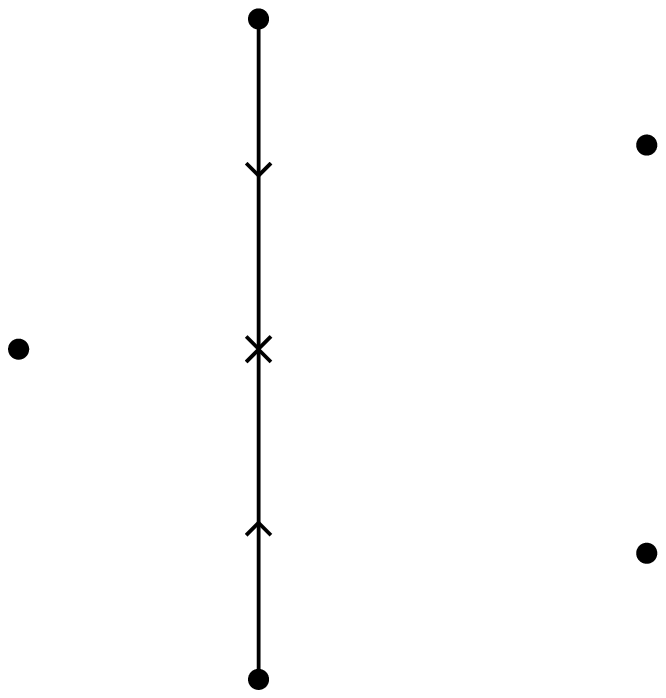}%
\\
\hfil $\QQQ^-(z)$ \hfil
}
\caption{The movement of the branching points along $\gamma^\pm_R$}
\label{fig3}
\end{figure}
\begin{figure}
\includegraphics[scale=.6]{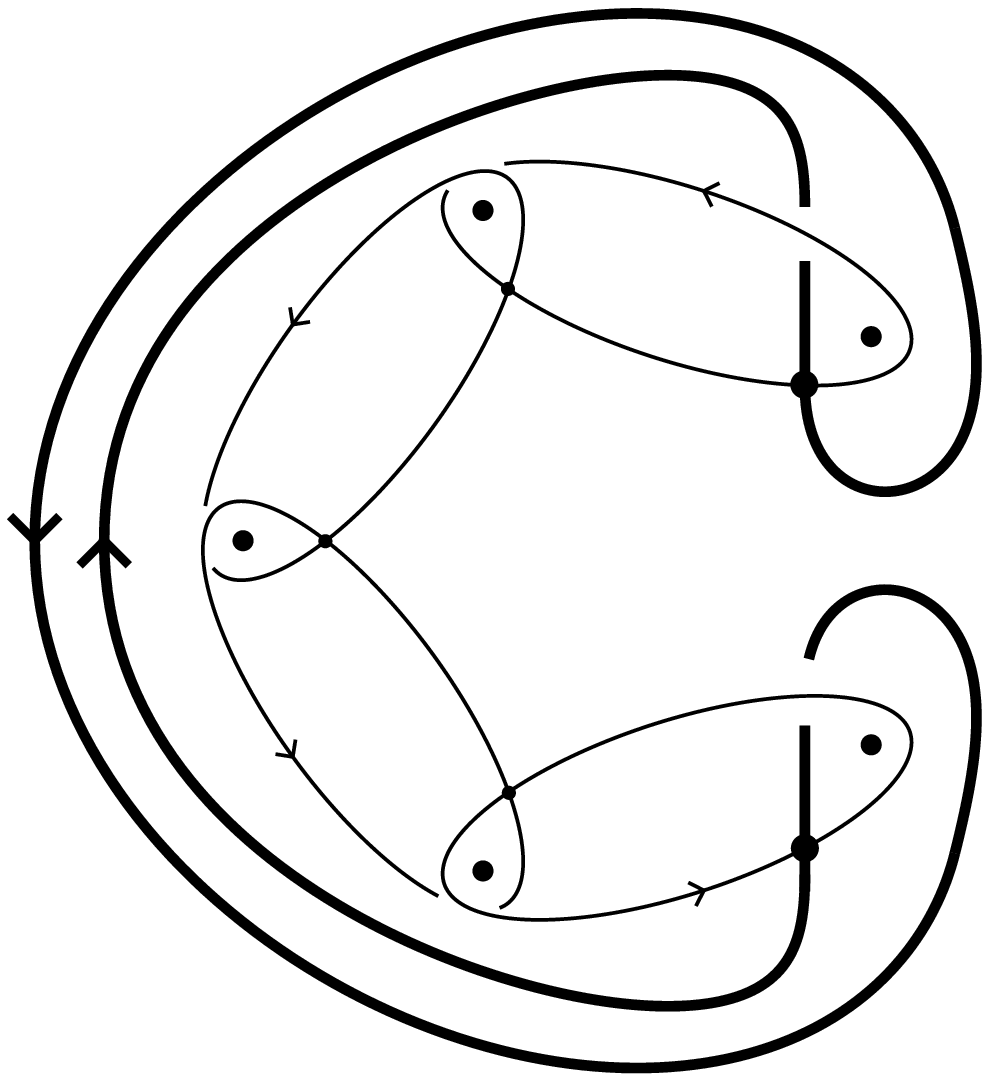}%
\caption{The vanishing cycle $\sigma_R^+=-\bdr \Gamma_R\sp+ $}
\label{figsigma}
\end{figure}
\begin{figure}
\parbox{170pt}{
\includegraphics[scale=.55]{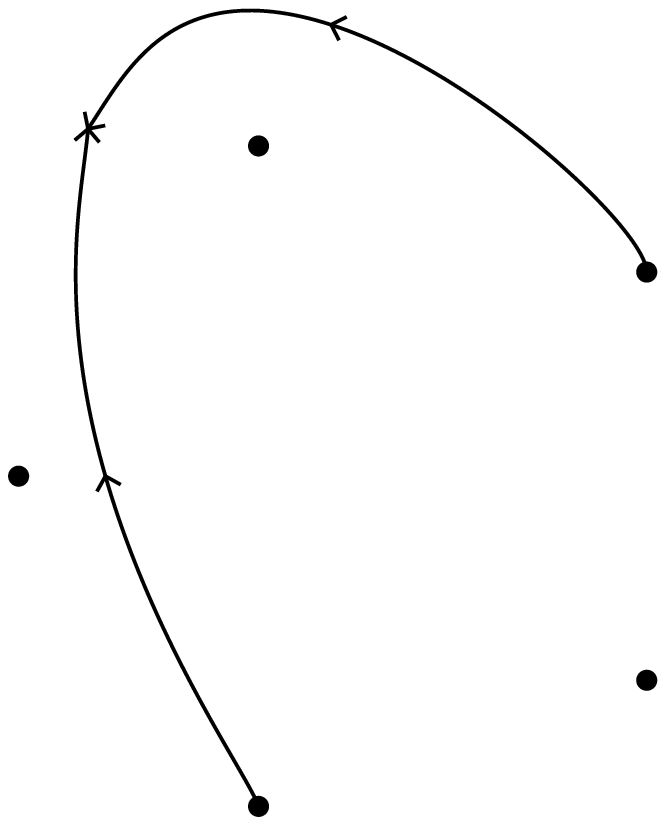}%
\\
\hfil $\QQQ^+(z)$ \hfil
}
\vrule height 100pt 
\parbox{170pt}{
\includegraphics[scale=.55]{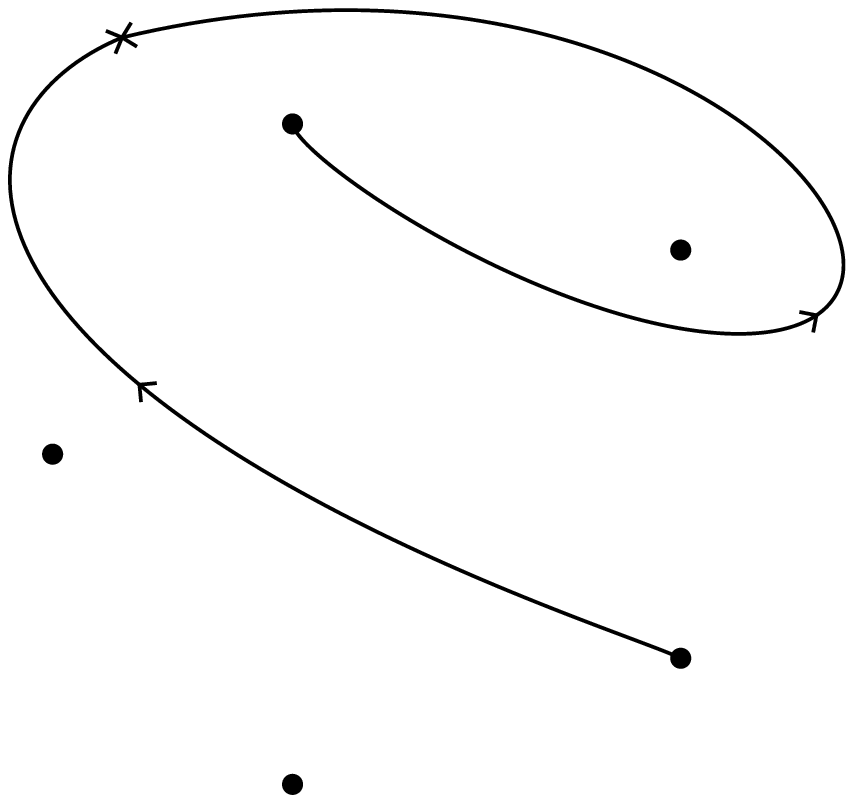}%
\\
\hfil $\QQQ^-(z)$ \hfil
}
\caption{The movement of the branching points along $\gamma^\pm_S$}
\label{fig4}
\end{figure}

\medskip
Since $Q\sp\pm$ is smooth at the origin and intersects $L$ with multiplicity $5$ at the origin, 
the movement of the branching points $\QQQ^\pm(z)$ along the loop $\lambda$ 
is homotopically equivalent to the rotation around the origin of the angle $2\pi/5$.
Hence the  monodromy on the simple closed curves 
is given by
$a_\nu\mapsto a_{\nu+1}$.
Let 
$$
\Lambda_\nu\;:\; \sphere^1\times I\to W^\pm
$$ 
be the topological  $2$-chain
over $\lambda$ that connects $a_\nu$ and  $a_{\nu+1}$.
We have
$$
\bdr [\Lambda_\nu]=[a_{\nu+1}]-[a_\nu].
$$
%
%
%
%
The movement of  the  branching points  $\QQQ^\pm (z)$
when $z$ moves from $b$ to $R^\pm$ along the path $\gamma_R^\pm$
is homotopically equivalent to the movement depicted in 
 Figure~\ref{fig3}.
 Let 
 $$
 \Gamma_R\sp\pm \;:\;\Theta\to W^\pm
 $$
  be the thimble over $\gamma_R^\pm$
corresponding to the critical point of $q^\pm: W^\pm\to U$ in the fiber over  $R^\pm$.
The vanishing cycle
$\sigma_R^+=-\bdr \Gamma_R\sp+ $
is depicted by a thick line in Figure~\ref{figsigma}.
We  choose the orientation
of $\sigma_R^+$ as in Figure~\ref{figsigma}.
Then we have 
$$
([\sigma_R^+], [a_1])_F=1,
\quad
([\sigma_R^+], [a_2])_F=([\sigma_R^+], [a_3])_F=0,
\quad
([\sigma_R^+], [a_4])_F=1,
$$ 
and hence
$$
[\sigma_R^+]=[a_1]-[a_2]+[a_3]-[a_4].
$$
In the same way, we see that the homology class of the vanishing cycle $\sigma_R^-=-\bdr \Gamma_R\sp-$ is equal to 
$$
[\sigma_R^-]=[a_2]+[a_3]
$$
under an appropriate choice of orientation.
The movement of  the points $\QQQ^\pm (z)$
when $z$ moves from $b$ to $S^\pm$ along the path $\gamma_S^\pm$
is homotopically equivalent to the movement depicted in 
 Figure~\ref{fig4}.
We choose the orientations 
of the thimbles 
$$
\Gamma_S^\pm\;:\; \Theta\to W^\pm 
$$
over $\gamma_S^\pm$  in such a way that  the homology classes of the vanishing cycles 
$\sigma_S^\pm=-\bdr \Gamma_S^\pm$ 
are
\begin{eqnarray*}
{} [\sigma_S^+]&=&[a_1]-[a_2]-[a_3] \;\;\;\;\quand \\
{} [\sigma_S^-]&=&2[a_1]-[a_2]-[a_3]-[a_4].
\end{eqnarray*}
The movement of  the  points $\QQQ^\pm (z)$
for the path $\gamma_{\bar{S}}^\pm$
is obtained from  Figure~\ref{fig4} by the  conjugation $\bar{\phantom{z}}:\C\to\C$ over $\R$.
We  choose the orientations 
of the thimbles $\Gamma_{\bar{S}}^\pm$ in such a way that
 $[\sigma_{\bar S}^\pm]=-\bdr [\Gamma_{\bar S}^\pm]$ are equal to 
\begin{eqnarray*}
{} [\sigma_{\bar{S}}^+]&=&-[a_2]-[a_3]+[a_4] \;\;\;\;\quand \\
{} [\sigma_{\bar{S}}^-]&=&-[a_1]-[a_2]-[a_3]+2 [a_4].
\end{eqnarray*}
\par
\medskip
%
Now we can calculate the kernel $\homH_2(W^\pm, \Z)$ of the homomorphism
$$
\bdr\;\;:\;\;\bigoplus_{\nu=1}^4\Z[\Lambda_\nu]\oplus \Z[\Gamma^\pm_R]\oplus \Z[\Gamma^\pm_S]\oplus  \Z[\Gamma^\pm_{\bar S}]
\;\to\; \bigoplus _{\nu=1}^4\Z[a_\nu].
$$
We see that  $\homH_2(W^+, \Z)$ is  a free $\Z$-module of rank $3$ generated by
\begin{eqnarray*}
S^+_1&:=&-[\Lambda_1]-[\Lambda_3]+[\Gamma_R^+],\\
S^+_2&:=&-6[\Lambda_1]-2[\Lambda_2]+2[\Lambda_3]+[\Lambda_4]+5[\Gamma_S^+],\\
S^+_3&:=&[\Lambda_1]+[\Lambda_2]+[\Lambda_3]-[\Gamma_S^+]+[\Gamma_{\bar S}^+],
\end{eqnarray*}
while  $\homH_2(W^-, \Z)$ is  a free $\Z$-module of rank $3$ generated by
%
%
\begin{eqnarray*}
S^-_1&:=&-4[\Lambda_1]-3[\Lambda_2]-2[\Lambda_3]+[\Gamma_R^-]+2[\Gamma_S^-],\\
S^-_2&:=&-11[\Lambda_1]-7[\Lambda_2]-3[\Lambda_3]+[\Lambda_4]+5[\Gamma_S^-],\\
S^-_3&:=&3[\Lambda_1]+3[\Lambda_2]+3[\Lambda_3]-[\Gamma_S^-]+[\Gamma_{\bar S}^-].
\end{eqnarray*}
\par
\begin{figure}
\hskip -160pt 
\hbox{\includegraphics[scale=.8]{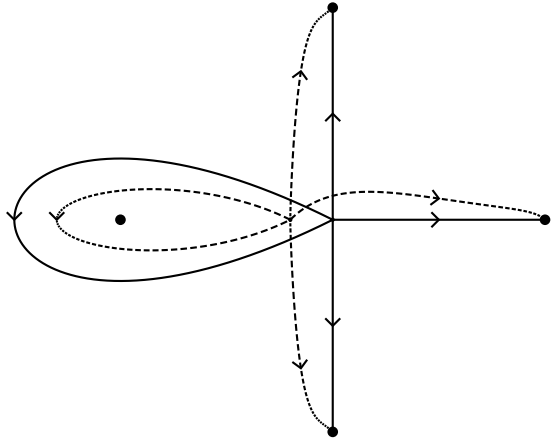}%
}
\caption{The perturbation}
\label{fig5}
\end{figure}

\medskip
We deform the loop $\lambda$ and 
the paths $\gamma^\pm_R$, $\gamma^\pm_S$ and $\gamma^\pm_{\bar{S}}$
as in Figure~\ref{fig5}.
The deformed loop $\lambda\sprime$ and paths  
$\gamma\sp{\prime\pm}_R$, $\gamma\sp{\prime\pm}_S$, $\gamma\sp{\prime\pm}_{\bar{S}}$
are depicted by the dotted curves.
Then the intersection numbers of the topological $2$-chains
$T=\Lambda_\nu, \Gamma^\pm_R, \Gamma^\pm_S, \Gamma^\pm_{\bar S}$ and  
$T\sprime=\Lambda\sprime_\nu, \Gamma^{\prime\pm}_R, \Gamma^{\prime\pm}_S, \Gamma^{\prime\pm}_{\bar S}$ 
are calculated 
as in Table~\ref{table:TT}.
Remark that 
the local intersection number $(t, t\sprime)_q$ 
of the underlying paths $t$ of $T$ and $t\sprime$ of $T\sprime$ is $-1$ for any 
intersection point $q$ contained in the small open neighborhood $\Delta$ of $b$ and $b\sprime$.
\begin{table}
$$
\renewcommand{\arraystretch}{1.7}
\begin{array}{c| cccc}
T\sprime \backslash T & \Lambda_\nu & \Gamma^\pm_R  & \Gamma^\pm_S  & \Gamma^\pm_{\bar S}\\
\hline
\Lambda\sprime_\nu  &0&0&0&0 \\
\Gamma^{\prime\pm}_R&([a_\nu], [\sigma_R\sp\pm])_F&-1&([\sigma_S\sp\pm], [\sigma_R\sp\pm])_F&0 \\
\Gamma^{\prime\pm}_S&([a_\nu], [\sigma_S\sp\pm])_F&0&-1&0 \\
\Gamma^{\prime\pm}_{\bar S} &([a_{\nu+1}], [\sigma_{\bar S}\sp\pm])_F&0&0&-1 \\
\end{array}
$$
\vskip 10pt 
\caption{The intersection numbers of $T$ and $T\sprime$}
\label{table:TT}
\end{table}
Therefore the intersection matrix of 
$\homH_2(W^+ , \Z)$
is calculated as follows:
$$
\renewcommand{\arraystretch}{1.2}
\begin{array}{c| ccc}
&S_1^+&S_2^+&S_3^+ \\
\hline
S_1^+ &0&0&0 \\
S_2^+ &0&40&-5 \\
S_3^+ &0&-5&2 
\end{array}.
$$
Then $I(W^+)$ is generated by $S_1^+$,
where
$I(W^+)\subset \homH_2(W^+)$ is the submodule
defined by~\eqref{eq:defI}.
Thus $\Tr(X^+)\cong \homH_2(W^+)/I(W^+)$ is generated by $S_2^+ + I(W^+) $ and $S_3^+ +I(W^+)$, 
and $\Tr(X^+)$ is isomorphic to
$$
L[40, -5, 2]\cong L[2, 1, 28].
$$
The intersection matrix of 
$\homH_2(W^- , \Z)$
is calculated as follows:
$$
\renewcommand{\arraystretch}{1.2}
\begin{array}{c| ccc}
&S_1^-&S_2^-&S_3^- \\
\hline
S_1^- &22&55&-22 \\
S_2^- &55&140&-55 \\
S_3^- &-22&-55&22
\end{array}.
$$
Then $I(W^-)$ is generated by $S_1^- - S_3^-$.
Therefore
$\Tr(X^-)\cong \homH_2(W^-)/I(W^-)$ is generated by $S_2^- +I(W^-) $ and $S_3^- + I(W^-)$, and 
 $\Tr(X^-)$ is isomorphic to
$$
L[140, -55, 22]\cong L[8, 3, 8].
$$
Thus Proposition~\ref{prop:theexample} is proved.
\begin{remark}
For the algorithm   to determine whether given two lattices of rank $2$ are
isomorphic or not,
see~\cite[Chapter 15]{MR1662447}.
\end{remark}
\section{The equations}\label{sec:equation}
In this section,  we construct  homogeneous  polynomials of degree $6$ defining  complex
projective plane curves that have singular  points of type $A_{10}$
and  of type $A_9$. 
Two of such polynomials are as follows:
\begin{eqnarray}
&10y^4z^2-20xy^2z^3+10x^2z^4-(-108\pm40\sqrt{5})\,c_0 xy^3z^2+&
\label{a10a9equ}\\
&(-64\pm20\sqrt{5}\,)c_0x^2yz^3+(-44\pm20\sqrt{5}\,)c_0y^5z-(-58\pm30\sqrt{5}\,){c_0}^2x^3z^3+&
\nonumber\\
&(-48\pm30\sqrt{5}\,){c_0}^2x^2y^2z^2+(-14\pm10\sqrt{5}\,){c_0}^3x^3yz^2+(-9\pm5\sqrt{5}\,){c_0}^4x^4z^2,&
\nonumber
\end{eqnarray}
with $c_0\in\mathbb{C}^\times$. We explain how to obtain these equations.
%
\par
\medskip
First we prove a lemma.
\begin{lemma}\label{A10}
Let $f(x, y)=0$ be a defining equation of a complex affine plane curve of  degree $6$ that has a
singular point of type $A_{10}$ at the origin with the tangent  $x=0$. 
Then, after appropriate coordinate change of the form $(x, y)\mapsto (x, a y)$
with $a\ne 0$, $f$ is equal to one of the following polynomials
~\eqref{a2345} or~\eqref{a345} up to multiplicative
constant:
\begin{eqnarray}\label{eq:A10_1}
\lefteqn{\phantom{f=}x^2-2xy^2+y^4+a_{0,5}y^5+a_{0,6}y^6+a_{1,3}xy^3+a_{1,4}xy^4}&&
\label{a2345}\\
&&+a_{1,5}xy^5+a_{2,1}x^2y+a_{2,2}x^2y^2+a_{2,3}x^2y^3+a_{2,4}x^2y^4+a_{3,0}x^3
\nonumber\\
&&+a_{3,1}x^3y+a_{3,2}x^3y^2+a_{3,3}x^3y^3+a_{4,0}x^4+a_{4,1}x^4y+a_{4,2}x^4y^2
\nonumber\\
&&+a_{5,0}x^5+a_{5,1}x^5y+a_{6,0}x^6, \nonumber
\end{eqnarray}
where
$$\begin{array}{ll}
a_{3,0}=&{c_0}^2-a_{0,6}-a_{1,4}-a_{2,2},\\
a_{2,1}=&-2c_0+a_{0,5},\\
a_{1,3}=&2c_0-2a_{0,5},\\
a_{4,0}=&-\frac{1}{\,2\,}(-3{c_0}^4+6{c_0}^2c_1-2{c_1}^2+3{c_0}^2a_{0,6}-4c_0a_{1,5}-{c_0}^2a_{2,2}\\
    &-2c_0a_{2,3}+2a_{2,4}+2a_{3,2}),\\
a_{3,1}=&-2{c_0}^3+2c_0c_1+{c_0}^2a_{0,5}-a_{1,5}-a_{2,3},\\
a_{1,4}=&-\frac{1}{\,2\,}({c_0}^2-2c_1-2c_0a_{0,5}+3a_{0,6}+a_{2,2}),\\
a_{5,0}=&-\frac{1}{\,2\,}(-16{c_0}^6+45{c_0}^4c_1-33{c_0}^2{c_1}^2+2{c_1}^3-12{c_0}^3c_2+16c_0c_1c_2-2{c_2}^2)\\
   
&-(2{c_0}^5-3{c_0}^3c_1+c_0{c_1}^2)a_{0,5}-\frac{1}{\,2\,}(8{c_0}^4-9{c_0}^2c_1+3{c_1}^2)a_{0,6}\\
   
&-(-3{c_0}^3+2c_0c_1)a_{1,5}-\frac{1}{\,2\,}({c_0}^2c_1-{c_1}^2)a_{2,2}-(3{c_0}^2-2c_1)a_{2,4}\\
    &-({c_0}^2-c_1)a_{3,2}+c_0a_{3,3}-a_{4,2},\\
a_{4,1}=&-\frac{1}{\,2\,}(15{c_0}^5-30{c_0}^3c_1+12c_0{c_1}^2+8{c_0}^2c_2-4c_1c_2)\\
   
&-(-4{c_0}^4+4{c_0}^2c_1-{c_1}^2)a_{0,5}-\frac{1}{\,2\,}{c_0}^3a_{0,6}-{c_0}^2a_{1,5}\\
    &+\frac{1}{\,2\,}{c_0}^3a_{2,2}+2c_0a_{2,4}+c_0a_{3,2}-a_{3,3},\\
a_{2,3}=&3{c_0}^3-4c_0c_1+2c_2-(3{c_0}^2-2c_1)a_{0,5}+3c_0a_{0,6}-2a_{1,5}-c_0a_{2,2},\\
a_{5,1}=&\frac{1}{\,2\,}(26{c_0}^7-83{c_0}^5c_1+74{c_0}^3{c_1}^2-12c_0{c_1}^3+25{c_0}^4c_2-40{c_0}^2c_1c_2\\
   
&+4{c_1}^2c_2+6c_0{c_2}^2-2c_3)-(2{c_0}^6-9{c_0}^4c_1+9{c_0}^2{c_1}^2+3{c_0}^3c_2-6c_0c_1c_2\\
   
&+{c_2}^2)a_{0,5}-\frac{1}{\,2}(-34{c_0}^5+59{c_0}^3c_1-24c_0{c_1}^2-9{c_0}^2c_2+6c_1c_2)a_{0,6}\\
   
&-(9{c_0}^4-11{c_0}^2c_1+2c_0c_2)a_{1,5}-\frac{1}{\,2\,}(2{c_0}^5-7{c_0}^3c_1+6c_0{c_1}^2+{c_0}^2c_2\\
   
&-2c_1c_2)a_{2,2}-(-7{c_0}^3+8c_0c_1-2c_2)a_{2,4}-(-2{c_0}^3+3c_0c_1-c_2)a_{3,2}\\
    &-(2{c_0}^2-c_1)a_{3,3}+c_0a_{4,2}, 
\end{array}$$
with $a_{i,j}$, $c_k\in\mathbb{C}$,  or
\begin{eqnarray}\label{A10_2}
\lefteqn{\phantom{f=}x^2-2xy^3+y^6-2x^2y{c_0}+2xy^4{c_0}+x^2y^2\big({c_0}^2+2({c_0}^2-{c_1})\big)}
\label{a345}\\
&&-2xy^5({c_0}^2-{c_1})+x^3\big(-2{c_0}({c_0}^2-{c_1})-a_{2,3}\big)+x^2y^3a_{2,3}\nonumber\\
&&+x^2y^4a_{2,4}+x^3y(3{c_0}^4-4{c_0}^2{c_1}+{c_1}^2+{c_0}a_{2,3}-a_{2,4})+x^3y^3a_{3,3}\nonumber\\
&&+x^3y^2\big(-3{c_0}^5+6{c_0}^3{c_1}-3{c_0}{c_1}^2+{c_2}-({c_0}^2-{c_1})a_{2,3}+{c_0}a_{2,4}\big)\nonumber\\
&&+x^4a_{4,0}+x^4ya_{4,1}+x^4y^2a_{4,2}+x^5a_{5,0}+x^5ya_{5,1}+x^6a_{6,0},\nonumber
\end{eqnarray}
with $a_{i,j}$, $c_k\in\mathbb{C}$.

Conversely,
if $c_3\ne 0$,
then the affine curve defined by the polynomial~\eqref{eq:A10_1} has a singular point of type $A_{10}$
at the origin,
and if $c_2\ne 0$,
then the affine curve defined by the polynomial~\eqref{A10_2} has a singular point of type $A_{10}$
at the origin.
\end{lemma}
We will use the following method to determine the type of singularities from
the form of the equation.
\begin{definition}
Let $k$ be an algebraically closed field and let
$w=(w_0,w_1)\in\mathbb{Q}_{\geq0}^{\ 2}$. Let $M=x^{e_0}y^{e_1}\in k[[x,y]]$
be a monomial. We define the {\it weight} of $M$ by $w(M):=\sum e_iw_i$. A
formal power series $f\in k[[x,y]]$ is said to be {\it semi-quasihomogeneous}
with respect to the weight $w$ if $f$ is of  the form
$f=f_{w=1}+f_{w>1}$ such that
\begin{list}{}{} 
\item[\textup{(i)}] every non-zero coefficient monomial $M$ in $f_{w=1}$
satisfies $w(M)=1$, and $f_{w=1}$ defines an isolated singularity, and 
\item[\textup{(ii)}] every non-zero coefficient monomial $M$ in $f_{w>1}$
satisfies $w(M)>1$. 
\end{list}
A semi-quasihomogeneous $f$ is said to be {\it quasihomogeneous} with respect to
the weight $w$ if $f_{w>1}=0$.
\end{definition}
\begin{proposition}[\cite{MR1175728}, Proposition 2.3]\label{recognition}
A semi-quasihomogeneous 
$f\in k[[x,y]]$ with respect to the weight ${\mathcal
A}_m=(\frac{1}{\,2\,},\frac{1}{m+1})$, ${\mathcal
D}_m=(\frac{1}{m-1},\frac{m-2}{2(m-1)})$, ${\mathcal
E}_6=(\frac{1}{\,3\,},\frac{1}{\,4\,})$, ${\mathcal
E}_7=(\frac{1}{\,3\,},\frac{2}{\,9\,})$ and ${\mathcal
E}_8=(\frac{1}{\,3\,},\frac{1}{\,5\,})$ defines a simple singularity if
$f_{w=1}$ defines an isolated singularity at the origin. The type of the singular point is
$A_m$, $D_m$, $E_6$, $E_7$ and $E_8$ respectively.
\end{proposition} 
\proof[Proof of Lemma \ref{A10}]{ Let $f_0(x,y)=\sum
b_{i,j}x^iy^j\in\mathbb{C}[x,y]$ be a polynomial of  degree $6$ with complex
coefficients $b_{i,j}$. Suppose that the affine plane curve defined by $f_0$ has a
singularity of type $A_{10}$ at $(0,0)$ with the tangent $x=0$. We can write 
$$f_0=x^2+b_{3,0}x^3+b_{2,1}x^2y+b_{1,2}xy^2+b_{0,3}y^3+({\rm higher\
terms}).$$

Firstly, let $w=(\frac{1}{\,2\,},\frac{1}{\,3\,})$. 
If $b_{0,3}\neq0$, then 
 $(f_0)_{w=1}=x^2+b_{0,3}y^3$  would define an isolated singularity at the origin and 
hence $f_0=0$ would have a singularity of type $A_2$ at the origin  by Proposition~\ref{recognition}. 
Thus $b_{0,3}$
must be equal to $0$. 

Secondly, let $w=(\frac{1}{\,2\,},\frac{1}{\,4\,})$.
Then  $f_0$ is
semi-quasihomogeneous with respect to  $w$, and hence the
quasihomogeneous part $(f_0)_{w=1}$ must define a \emph{non}-isolated singularity
at the origin by Proposition \ref{recognition}. Hence there exists
$b\in\mathbb{C}$ such that $(f_0)_{w=1}$ is equal to $x^2-2bxy^2+b^2y^4$. We
divide into two cases, the case where $b\neq0$ and the case where $b=0$.

{\it Case 1} $(b\neq0)$.  
We  change the coordinate via $\sqrt{b}\,y\mapsto y$.
 Then  we have
$(f_0)_{w=1}\mapsto x^2-2xy^2+y^4$. 
Therefore,  without loss of generality, we can write
\begin{eqnarray}
\lefteqn{f_0(x,y)=x^2-2xy^2+y^4+a_{0,5}y^5+a_{0,6}y^6+a_{1,3}xy^3+a_{1,4}xy^4}&&
\label{standardform}\\
&&+a_{1,5}xy^5+a_{2,1}x^2y+a_{2,2}x^2y^2+a_{2,3}x^2y^3+a_{2,4}x^2y^4+a_{3,0}x^3
\nonumber\\
&&+a_{3,1}x^3y+a_{3,2}x^3y^2+a_{3,3}x^3y^3+a_{4,0}x^4+a_{4,1}x^4y+a_{4,2}x^4y^2
\nonumber\\
&&+a_{5,0}x^5+a_{5,1}x^5y+a_{6,0}x^6, \nonumber
\end{eqnarray}
with $a_{i,j}\in\mathbb{C}$. 

Change the coordinate via $x\mapsto x+y^2$. Assume that this coordinate change
transforms $f_0$ into $f_1$. An elementary calculation shows that
\begin{eqnarray*}
\lefteqn{f_0\mapsto
f_1=x^2+a_{2,1}x^2y+(a_{0,5}+a_{1,3}+a_{2,1})y^5+(a_{1,3}+2a_{2,1})xy^3+a_{3,0}x^3}&&
\\
&&+(a_{0,6}+a_{1,4}+a_{2,2}+a_{3,0})y^6+(a_{2,2}+3a_{3,0})x^2y^2+(a_{1,4}+2a_{2,2}+3a_{3,0})xy^4\\
&&+a_{3,1}x^3y+(a_{1,5}+a_{2,3}+a_{3,1})y^7+(a_{2,3}+3a_{3,1})x^2y^3+(a_{1,5}+2a_{2,3}+3a_{3,1})xy^5\\
&&+a_{4,0}x^4+(a_{2,4}+a_{3,2}+a_{4,0})y^8+(a_{3,2}+4a_{4,0})x^3y^2+(2a_{2,4}+3a_{3,2}+4a_{4,0})xy^6\\
&&+(a_{2,4}+3a_{3,2}+6a_{4,0})x^2y^4+a_{4,1}x^4y+(a_{3,3}+a_{4,1})y^9+3(a_{3,3}+2a_{4,1})x^2y^5\\
&&+(a_{3,3}+4a_{4,1})x^3y^3+(3a_{3,3}+4a_{4,1})xy^7+a_{5,0}x^5+(a_{4,2}+a_{5,0})y^{10}\\
&&+(a_{4,2}+5a_{5,0})x^4y^2+2(2a_{4,2}+5a_{5,0})x^3y^4+2(3a_{4,2}+5a_{5,0})x^2y^6\\
&&+(4a_{4,2}+5a_{5,0})xy^8+a_{5,1}x^5y+5a_{5,1}x^4y^3+10a_{5,1}x^3y^5+10a_{5,1}x^2y^7\\
&&+5a_{5,1}xy^9+a_{5,1}y^{11}+a_{6,0}x^6+6a_{6,0}x^5y^2+15a_{6,0}x^4y^4+20a_{6,0}x^3y^6\\
&&+15a_{6,0}x^2y^8+6a_{6,0}xy^{10}+a_{6,0}y^{12}.
\end{eqnarray*}
First, 
let $w=(\frac{\,\,\,1}{\,2\,},\frac{1}{\,5\,})$. If there were the term $y^5$
with non-zero coefficient in $f_1$, the singularity of $f_1=0$ at the origin would be of type $A_4$
by Proposition \ref{recognition}. 
Next, let
$w=(\frac{1}{\,2\,},\frac{1}{\,6\,})$. By the same argument as above, the
quasihomogeneous part $(f_1)_{w=1}$ of $f_1$ must define non-isolated
singularities at the origin, because  otherwise $f_1=0$ would have a singularity of type
$A_5$. Thus there exists a complex number $c_0$ such that
$(f_1)_{w=1}=x^2+(a_{1,3}+2a_{2,1})xy^3+(a_{0,6}+a_{1,4}+a_{2,2}+a_{3,0})y^6$
is equal to    $x^2-2c_0xy^3+{c_0}^2y^6$. Consequently we have following
conditions:
\begin{list}{}{}
\item [{\it Step 1.}] ($x\mapsto x+y^2$)\\
$a_{0,5}+a_{1,3}+a_{2,1}=0$, \\
$a_{0,6}+a_{1,4}+a_{2,2}+a_{3,0}={c_0}^2$ and \\
$a_{1,3}+2a_{2,1}=-2c_0$. 
\end{list}

Then we change the coordinate via $x\mapsto x+c_0y^3$. Assume that this
transformation takes $f_1$ to $f_2$. Let $w=(\frac{1}{\,2\,},\frac{1}{\,7\,})$.
The coefficient of $y^7$ is equal to $0$, because  otherwise the singularity at the
origin would be of type $A_{6}$. Next, let
$w=(\frac{1}{\,2\,},\frac{1}{\,8\,})$. The quasihomogeneous part $(f_2)_{w=1}$
must define non-isolated singularities at the origin. 
Hence there exists $c_1\in\mathbb{C}$ such
that $(f_2)_{w=1}$ is equal to $x^2-2c_1xy^4+{c_1}^2y^8$. Therefore we have
\begin{list}{}{}
\item [{\it Step 2.}] ($x\mapsto x+c_0y^3$)\\
${c_0}^3+{c_0}^2a_{0,5}-3c_0a_{0,6}-2c_0a_{1,4}+a_{1,5}-c_0a_{2,2}+a_{2,3}+a_{3,1}=0$,
\\
$3{c_0}^4-3{c_0}^2a_{0,6}-3{c_0}^2a_{1,4}+c_0a_{1,5}-2{c_0}^2a_{2,2}+2c_0a_{2,3}+a_{2,4}+3c_0a_{3,1}+a_{3,2}
+a_{4,0}={c_1}^2$ and \\
$-{c_0}^2+2c_0a_{0,5}-3a_{0,6}-2a_{1,4}-a_{2,2}=-2c_1$.
\end{list}

The coordinate change via $x\mapsto x+c_1y^4$ takes
$f_2$ to $f_3$. The coefficient of $y^9$ in $f_3$ is equal to $0$ and there
exists $c_2\in\mathbb{C}$ such that $(f_3)_{w=1}=x^2-2c_2xy^5+{c_2}^2y^{10}$, 
where  $w=(\frac{1}{\,2\,},\frac{1}{\,10\,})$.
\begin{list}{}{}
\item [{\it Step 3.}] ($x\mapsto x+c_1y^4$)\\
$\frac{\,3\,}{2}{c_0}^5-4c_1{c_0}^3+2{c_0}^4a_{0,5}-3c_1{c_0}^2a_{0,5}+{c_1}^2a_{0,5}-\frac{\,11\,}{2}{c_0}^3a_{0,6}+3c_1c_0a_{0,6}+5{c_0}^2a_{1,5}+\frac{\,3\,}{2}{c_0}^3a_{2,2}-a_{2,2}c_1c_0+2{c_0}^2a_{2,3}-a_{2,3}c_1-2c_0a_{2,4}-c_0a_{3,2}+a_{3,3}+a_{4,1}=0$,
\\
$7{c_0}^6-\frac{\,35\,}{2}c_1{c_0}^4+\frac{\,15\,}{2}{c_1}^2{c_0}^2+{c_1}^3+{c_0}^5a_{0,5}+3c_1{c_0}^3a_{0,5}-3{c_1}^2c_0a_{0,5}-9{c_0}^4a_{0,6}-\frac{\,9\,}{2}c_1{c_0}^2a_{0,6}+\frac{\,3\,}{2}{c_1}^2a_{0,6}+11{c_0}^3a_{1,5}+2c_1c_0a_{1,5}+3{c_0}^4a_{2,2}+\frac{1}{\,2\,}c_1{c_0}^2a_{2,2}-\frac{1}{\,2\,}{c_1}^2a_{2,2}+5{c_0}^3a_{2,3}-5{c_0}^2a_{2,4}-2c_1a_{2,4}-3{c_0}^2a_{3,2}-c_1a_{3,2}+3c_0a_{3,3}+4c_0a_{4,1}+a_{4,2}+a_{5,0}={c_2}^2$
and \\
$3{c_0}^3-4c_1c_0-3{c_0}^2a_{0,5}+2c_1a_{0,5}+3c_0a_{0,6}-2a_{1,5}-c_0a_{2,2}-a_{2,3}=-2c_2$.
\end{list}

Change the coordinate via $x\mapsto x+c_2y^5$. Suppose that this transformation
takes $f_3$ to $f_4$. Let $w=(\frac{1}{\,2\,},\frac{1}{\,11\,})$. 
Then   $(f_4)_{w=1}$ is
equal to $x^2+c_3y^{11}$,
where $c_3$ is given below.
If $c_3\ne 0$,
then $(f_4)_{w=1}=0$ defines an isolated singular point at the origin,
and hence $f_4=0$ has a singular point of type $A_{10}$ at the origin.
\begin{list}{}{}
\item [{\it Step 4.}] ($x\mapsto x+c_2y^5$)\\
$13{c_0}^7-\frac{\,83\,}{2}c_1{c_0}^5+37{c_1}^2{c_0}^3-6{c_1}^3c_0+\frac{\,25\,}{2}c_2{c_0}^4-20c_2c_1{c_0}^2+c_2{c_1}^2+3{c_2}^2c_0+2{c_0}^6a_{0,5}-9c_1{c_0}^4a_{0,5}+9{c_1}^2{c_0}^2a_{0,5}+3c_2{c_0}^3a_{0,5}-6c_2c_1c_0a_{0,5}+{c_2}^2a_{0,5}-17{c_0}^5a_{0,6}+\frac{\,59\,}{2}c_1{c_0}^3a_{0,6}-12{c_1}^2c_0a_{0,6}-\frac{\,9\,}{2}c_2{c_0}^2a_{0,6}+3c_2c_1a_{0,6}+9{c_0}^4a_{1,5}-11c_1{c_0}^2a_{1,5}+{c_1}^2a_{1,5}+2c_2c_0a_{1,5}+{c_0}^5a_{2,2}-\frac{\,7\,}{2}c_1{c_0}^3a_{2,2}+3{c_1}^2c_0a_{2,2}+\frac{1}{\,2\,}c_2{c_0}^2a_{2,2}-c_2c_1a_{2,2}-7{c_0}^3a_{2,4}+8c_1c_0a_{2,4}-2c_2a_{2,4}-2{c_0}^3a_{3,2}+3c_1c_0a_{3,2}-c_2a_{3,2}+2{c_0}^2a_{3,3}-c_1a_{3,3}-c_0a_{4,2}+a_{5,1}=c_3$.
\end{list}

Solve the system of linear equations appearing in each step by choosing
unknowns suitably. Then we have the solutions denoted in (\ref{a2345}).

{\it Case 2} $(b=0)$. In this case, without loss of generality we can write
\begin{eqnarray*} 
\lefteqn{f_0=x^2-2b\sprime xy^3+{b\sprime }^2y^6+a_{1,4}xy^4+a_{1,5}xy^5+a_{2,1}x^2y+a_{2,2}x^2y^2}\\
&&+a_{2,3}x^2y^3+a_{2,4}x^2y^4+a_{3,0}x^3+a_{3,1}x^3y+a_{3,2}x^3y^2+a_{3,3}x^3y^3\\
&&+a_{4,0}x^4+a_{4,1}x^4y+a_{4,2}x^4y^2+a_{5,0}x^5+a_{5,1}x^5y+a_{6,0}x^6.
\end{eqnarray*}
Assume that $b\sprime =0$. If $a_{1,4}\neq0$, then the polynomial $f_0$ is
semi-quasihomogeneous with respect to the weight
$w=(\frac{1}{\,2\,},\frac{1}{\,8\,})$ and $(f_0)_{w=1}$ defines an isolated
singularity. Hence $f_0=0$ would have a singularity of type $A_7$ at the origin.
If $a_{1,4}=0$ and $a_{1,5}\neq0$, then $f_0$
is semi-quasihomogeneous with respect to 
$w=(\frac{1}{\,2\,},\frac{1}{\,10\,})$, and $(f_0)_{w=1}$ defines an isolated
singularity at the origin, so that  $f_0=0$ would have a
singularity of type $A_9$ at the origin. If
$a_{1,4}=a_{1,5}=0$, then $f_0$ defines non-isolated singularities at the origin. Therefore
$b\sprime $ is not equal to zero. Furthermore, the coordinate change
$\sqrt[3]{b\sprime }\,y\mapsto y$ takes $b\sprime $ to $1$. Therefore we can write
\begin{eqnarray}
\lefteqn{f_0=x^2-2xy^3+y^6+a_{1,4}xy^4+a_{1,5}xy^5+a_{2,1}x^2y+a_{2,2}x^2y^2+a_{2,3}x^2y^3}&&\label{eqn345}\\
&&+a_{2,4}x^2y^4+a_{3,0}x^3+a_{3,1}x^3y+a_{3,2}x^3y^2+a_{3,3}x^3y^3+a_{4,0}x^4
\nonumber\\
&&+a_{4,1}x^4y+a_{4,2}x^4y^2+a_{5,0}x^5+a_{5,1}x^5y+a_{6,0}x^6.\nonumber
\end{eqnarray} 
By a similar argument as in Case 1, we have the following three steps:
\begin{list}{}{}
\item [{\it Step 1.}] ($x\mapsto x+y^3$)\\
$a_{1,4}+a_{2,1}=0$,\\
$a_{1,5}+a_{2,2}+a_{3,0}={c_0}^2$ and \\
$a_{1,4}+2a_{2,1}=-2c_0$,
\item [{\it Step 2.}] ($x\mapsto x+c_0y^4$)\\
$-c_0a_{1,5}+a_{2,3}+a_{3,0}=0$, \\
${c_0}^4-{c_0}^2a_{1,5}+2c_0a_{2,3}+a_{2,4}+3c_0a_{3,0}+a_{3,1}={c_1}^2$ and
\\
$-2{c_0}^2-a_{1,5}=-2c_1$,
\item [{\it Step 3.}] ($x\mapsto x+c_1y^5$)\\
$3{c_0}^5-6{c_0}^3c_1+3c_0{c_1}^2+{c_0}^2a_{2,3}-c_1a_{2,3}-c_0a_{2,4}+a_{3,2}=c_2$,
\end{list}
where $c_i\in\mathbb{C}$. Regard $a_{1,4}$, $a_{1,5}$, $a_{2,1}$, $a_{2,2}$,
$a_{3,0}$, $a_{3,1}$ and $a_{3,2}$ as unknowns and solve the system of linear
equations. The solutions are
\begin{eqnarray*}
a_{2,1}&=&-2c_0,\\
a_{2,2}&=&{c_0}^2-a_{1,5},\\
a_{1,4}&=&2c_0,\\
a_{1,5}&=&-2({c_0}^2-{c_1}),\\
a_{3,0}&=&-2{c_0}({c_0}^2-{c_1})-a_{2,3},\\
a_{3,1}&=&3{c_0}^4-4{c_0}^2{c_1}+{c_1}^2+{c_0}a_{2,3}-a_{2,4}\ {\rm and}\\
a_{3,2}&=&-3{c_0}^5+6{c_0}^3{c_1}-3{c_0}{c_1}^2+{c_2}-({c_0}^2-{c_1})a_{2,3}+{c_0}a_{2,4}.
\end{eqnarray*}
Substituting them for the coefficient of (\ref{eqn345}),  we obtain the
polynomial (\ref{a345}). 
\qed}\\

\begin{claim}
Let $F(x,y,z)\in\mathbb{C}[x,y,z]$ be a  homogeneous polynomial of degree $6$ that 
satisfies
$$F(x,y,1)=f(x,y),$$
where $f$ is the polynomial (\ref{a2345}) in the statement of Lemma \ref{A10} with $c_3\ne 0$.
Let $g(y,z):=F(1,y,z)$. Then $g$ is semi-quasihomogeneous with respect to the
weight $w=(\frac{1}{\,10\,},\frac{1}{\,2\,})$ if and only if 
\begin{eqnarray}\label{a0}
a_{0,6}=a_{1,5}=a_{2,4}=a_{3,2}=a_{3,3}=a_{4,2}=a_{6,0}=0
\end{eqnarray}
and
\begin{eqnarray}
\left\{\begin{array}{rcl} \label{c0}
c_1&=&\frac{1}{\,2\,}(5\pm\sqrt{5}\,){c_0}^2,\\
a_{2,2}&=&2(c_1+c_0a_{0,5})-{c_0}^2,\\
a_{0,5}&=&\frac{2}{\,5\,}(-11\pm5\sqrt{5}\,)c_0,\\
c_3&=&-\frac{6}{\,25\,}{(-123\pm55\sqrt{5}\,)}{c_0}^7.
\end{array}\right.
\end{eqnarray}
Moreover, if (\ref{a0}) and (\ref{c0}) hold, then $g_{w=1}$ defines an
isolated singularity at the origin
and hence $g=0$ has a singular point of type $A_9$ at $(0,0)$
by Proposition \ref{recognition}.
\end{claim}

\proof{
We write  $g$ in
the form
$$g=g_{w<1}+g_{w=1}+g_{w>1},$$
where $w=(\frac{\,1\,}{\,10\,}, \frac{\,1\,}{\,2\,})$.
The condition
$g_{w<1}=0$ is equivalent to (\ref{a0}) and
$$\left\{
\begin{array}{rcl}
0&=&2(c_1+c_0a_{0,5})-{c_0}^2-a_{2,2}\\
0&=&-4{c_0}^3+6c_1c_0+5a_{0,5}{c_0}^2-2a_{0,5}c_1-2c_2\\
0&=&5{c_0}^4-5{c_0}^2c_1+{c_1}^2\\
0&=&{a_{0,5}}^2(5{c_0}^2-2c_1)(15{c_0}^3-10c_0c_1+5{c_0}^2a_{0,5}-2c_1a_{0,5})-4c_3\\
0&=&12{c_0}^5-20{c_0}^3c_1+8c_0{c_1}^2+25{c_0}^4a_{0,5}-20{c_0}^2c_1a_{0,5}+4{c_1}^2a_{0,5}.
\end{array}\right.
$$
Solving  this system of equations,  we get~\eqref{c0}. Note that
we have $c_0\neq0$ by the assumption $c_3\neq0$.
Substituting \eqref{a0} and  \eqref{c0} for coefficients of $g$, we have 
$$g_{w=1}=\frac{-9\pm5\sqrt{5}\,}{10}{c_0}^4\,z^2+\frac{2(-11\pm5\sqrt{5}\,)}{5}c_0\,zy^5.$$
Since  $c_0\neq0$, $g_{w=1}$  defines an isolated
singularity. 
\qed}
\\

Note that $10F(x,y,z)$ is equal to (\ref{a10a9equ}) under the condition
(\ref{a0}) and (\ref{c0}). 
\par
\medskip
Finally, let $c_0=1$. The curve defined by the equation (\ref{a10a9equ}) has a
singular point of type $A_{10}$ at $(0:0:1)$, a singular point of type $A_9$ at
$(1:0:0)$,  and is smooth except for these two points.
\begin{remark}
In~\cite{MR1900779},
a different method to obtain  defining equations of sextic curves with big Milnor number
is given.
\end{remark}
\bibliographystyle{plain}
\def\cprime{$'$} \def\cprime{$'$} \def\cprime{$'$} \def\cprime{$'$}

\end{document}